\documentclass[11pt]{article}

\usepackage{amsmath,amssymb,amsfonts}
\usepackage{pstricks}

\numberwithin{equation}{section}
\newtheorem{thm}{Theorem}[section] 
\newtheorem{prp}[thm]{Proposition}
\newtheorem{lmm}[thm]{Lemma}   
\newtheorem{ques}[thm]{Question} 
 
\def\Bbb#1{\mathbb{#1}}
\def\frak{\mathfrak}
\def\e_ref#1{(\ref{#1})}
\def\under#1{\underline{#1}}
\def\ov#1{\overline{#1}}

\def\lra{\longrightarrow}
\def\Llra{\Longleftrightarrow}

\def\al{\alpha}
\def\be{\beta}
\def\ga{\gamma}
\def\la{\lambda}
\def\om{\omega}
\def\si{\sigma}
\def\vph{\varphi}

\def\Ga{\Gamma}
\def\La{\Lambda}
\def\Si{\Sigma}

\def\lan{\langle}
\def\ran{\rangle}
\def\lr#1{\lan{#1}\ran}
\def\blr#1{\big\lan{#1}\big\ran}

\def\i{\infty}
\def\eset{\emptyset}
\def\PP{\Bbb{P}^2}
\def\P{\Bbb{P}^n}
\def\R{\mathbb R}
\def\Hom{\textnormal{Hom}}

\begin{document}

\title{Counting Plane Rational Curves:\linebreak 
Old and New Approaches}
\author{Aleksey Zinger
\thanks{Supported by an NSF Postdoctoral Fellowship}}
\date{\today}
\maketitle

\begin{abstract}
\noindent
These notes are intended as an easy-to-read supplement to part of 
the background material presented in my talks on enumerative geometry.
In particular, the numbers $n_3$ and $n_4$ of plane rational cubics 
through eight points and of plane rational quartics through eleven points
are determined via the classical approach of counting curves.
The computation of the latter number also illustrates my topological approach 
to counting the zeros of a fixed vector bundle section that lie in the main stratum 
of a compact space.
The arguments used in the computation of the number $n_4$ extend easily to 
counting plane curves with two or three nodes, for example.
Finally, an inductive formula for the number $n_d$ of
plane degree-$d$ rational curves passing through $3d\!-\!1$ points
is derived via the modern approach of counting stable maps.
This method is far simpler.
\end{abstract}

\thispagestyle{empty}
\tableofcontents
\listoftables

\section{Introduction}
\label{intro_sec}

\noindent
Enumerative geometry of algebraic varieties is a field of 
mathematics that dates back to the nineteenth century.
The general goal of this subject is to determine 
the number of geometric objects that satisfy 
pre-specified geometric conditions.
The objects are typically (complex) curves in 
a smooth algebraic manifold.
Such curves are usually required to represent the given homology class, 
to have certain singularities, and 
to satisfy various contact conditions with respect to
a collection of subvarieties.
One of the most well-known examples of an enumerative problem~is

\begin{ques}
\label{g0n2_ques}
If $d$ is a positive integer, what is the number $n_d$ of degree-$d$ 
rational curves that pass through $3d\!-\!1$ points in general position
in the complex projective plane~$\PP$?
\end{ques}

\noindent
Since the number of (complex) lines through any two distinct points
is one, $n_1\!=\!1$.
A little bit of algebraic geometry and topology gives
$n_2\!=\!1$ and $n_3\!=\!12$; see Section~\ref{low_deg_sec}.
It is far harder to find that $n_4\!=\!620$, but 
this number was computed as early as the middle of the nineteenth century;
\hbox{see~\cite[p378]{Ze}}.
We give a ``classical-style" computation of this number in Section~\ref{deg4_sec}.
Along the way, we determine the number of plane quartics that pass through 12~points 
and have two nodes and the number of plane quartics that pass through 11~points 
and have a cusp and a simple node; see Table~\ref{quartics_table}.
The derivations of Subsections~\ref{p1_subs}-\ref{p3_subs}
easily extend to counting arbitrary-degree plane curves with two nodes, a node and a cusp,
and with three nodes; see Table~\ref{curvegen_table} for explicit formulas.
These curves are of course not rational in general.
Subsections~\ref{p2_subs} and~\ref{p3_subs} also
illustrate our approach to determining the number of zeros
of a fixed vector bundle section that lie in the main stratum of a space.
This approach is one of the two main tools that we have applied to 
a number of enumerative problems;
see~\cite{Z1} and~\cite{Z2}, for example.\\

\noindent
The higher-degree numbers $n_d$ remained unknown until the early~1990s,
when a recursive formula for the numbers~$n_d$ was announced 
in~\cite{KoMa} and~\cite{RuT}:
\begin{equation}\label{recursion_e_main}
n_d=\frac{1}{6(d\!-\!1)}\sum_{d_1+d_2=d}
\Bigg(\!d_1d_2\!-2\frac{(d_1\!-\!d_2)^2}{3d-2}\!\Bigg)
\binom{3d\!-\!2}{\!3d_1\!-\!1\!}\! d_1d_2n_{d_1}n_{d_2}.
\end{equation}
We describe the argument of the latter paper in Section~\ref{all_deg_sec}.
It can also be used to solve the natural generalization of Question~\ref{g0n2_ques} 
to the higher-dimensional projective spaces; see Section~10 in~\cite{RuT}.\\

\noindent
{\it Remark:} A derivation of~\e_ref{recursion_e_main}, which is classical in spirit,
appears in~\cite{Ra2} and is based on~\cite{Ra1}.
The approach of Section~\ref{deg4_sec} is more direct and involves no blowups.\\

\noindent
Subsection~\ref{n3_subs} and Section~\ref{deg4_sec}, 
which are not used in Section~\ref{all_deg_sec}, 
assume some familiarity with cohomology groups and chern classes.
All other non-elementary terms, including those used in Question~\ref{g0n2_ques}, 
are described in Appendix~\ref{basics_sec}.
A different (and far more extensive) introduction to enumerative
geometry, as well as to its relations with physics, is given in~\cite{Ka}.

\section{The Low-Degree Numbers}
\label{low_deg_sec}

\subsection{The Degree-One Number}
\label{n1_subs}

\noindent
We start by computing the number $n_1$ topologically.
Throughout these notes, we will use the homogeneous coordinates
$[X,Y,Z]$ on the complex projective plane of Question~\ref{g0n2_ques},
i.e.~we take
$$\PP=\big\{(X,Y,Z)\!\in\!\Bbb{C}^3\!: 
(X,Y,Z)\!\neq\!(0,0,0)\big\}\big/\Bbb{C}^*
=\big\{[X,Y,Z]\!: (X,Y,Z)\!\in\!\Bbb{C}^3\!-\!(0,0,0)\big\}.$$
In this section, we use the following lemma.

\begin{lmm}
\label{genus_lmm}
If $\ga\!\lra\!\PP$ is the tautological line bundle,
$d$ is positive integer, and  
$s\!\in\!\Ga(\PP;\ga^{*\otimes d})$ is transverse to the zero set,
the set $s^{-1}(0)$ is a smooth two-dimensional submanifold
of~$\PP$ of genus
$$g\big(s^{-1}(0)\big)=\binom{d\!-\!1}{2}.$$
\end{lmm}

\noindent
This lemma is proved in Subsection~\ref{taut_subs}.
It can easily be verified directly in the $d\!=\!1$ and $d\!=\!2$ cases.\\

\noindent
A line, or degree-one curve, in $\PP$ is the quotient by
the $\Bbb{C}^*$-action of the zero set of 
a nonzero homogeneous polynomial 
$$s_{a_{100},a_{010},a_{001}}\equiv a_{100}X+a_{010}Y+a_{001}Z$$ 
of degree one on~$\Bbb{C}^3\!-\!\{0\}$. 
In other words, a degree-one curve in $\PP$ has the~form
$${\cal C}={\cal C}_{a_{100},a_{010},a_{001}}=
\big\{[X,Y,Z]\!\in\!\PP\!: 
a_{100}X\!+\!a_{010}Y\!+\!a_{001}Z\!=\!0\big\}$$
for some $\big(a_{100},a_{010},a_{001}\big)\!\in\!\Bbb{C}^3\!-\!\{0\}$.
Furthermore, 
$${\cal C}_{a_{100},a_{010},a_{001}}=
{\cal C}_{b_{100},b_{010},b_{001}} \quad\Llra\quad
\big(a_{100},a_{010},a_{001}\big)=
\la\big(b_{100},b_{010},b_{001}\big)
~~\hbox{for some}~~ \la\!\in\!\Bbb{C}^*.$$
Thus, the space of all degree-one curves in $\PP$ is 
$${\cal D}_1=\big\{(a_{100},a_{010},a_{001})\!:
(a_{100},a_{010},a_{001})\!\neq\!(0,0,0)\big\}\big/\Bbb{C}^*\approx\PP.$$
\\

\noindent
A homogeneous polynomial $s\!=\!a_{100}X\!+\!a_{010}Y\!+\!a_{001}Z$ 
of degree one on $\Bbb{C}^3$ determines a section 
$s_{a_{100},a_{010},a_{001}}$ of the bundle $\ga^*\!\lra\!\PP$. 
If $(a_{100},a_{010},a_{001})\!\neq\!(0,0,0)$, this section
is transverse to the zero set.
Thus, by Lemma~\ref{genus_lmm}, for all 
$[a_{100},a_{010},a_{001}]\!\in\!{\cal D}_1$
the genus of ${\cal C}_{a_{100},a_{010},a_{001}}$ is zero,
i.e.~this is a rational curve.\\

\noindent
Finally, let $p_1\!=\![X_1,Y_1,Z_1]$ and $p_2\!=\![X_2,Y_2,Z_2]$ 
be two distinct points in~$\PP$.
The curve ${\cal C}_{a_{100},a_{010},a_{001}}$ passes through
the point~$p_i$ if and only if $s_{a_{100},a_{010},a_{001}}(p_i)\!=\!0$.
Thus, the number~$n_1$ is the number of elements
$[a_{100},a_{010},a_{001}]\!\in\!{\cal D}_1$ such that
\begin{equation}\label{n1_e}
\begin{cases}
a_{100}X_1+a_{010}Y_1+a_{001}Z_1=0;\\
a_{100}X_2+a_{010}Y_2+a_{001}Z_2=0.
\end{cases}
\end{equation}
The solution of each of these equations on ${\cal D}_1$ is a line.
Since $[X_1,Y_1,Z_1]\!\neq\![X_2,Y_2,Z_2]$, the two lines are distinct.
Since two lines in a plane, or~$\PP$, intersect in a single point,
$n_1\!=\!1$.
Stated differently, $n_1\!=\!1$ because the space of solutions
of the system~\e_ref{n1_e} in 
$(a_{100},a_{010},a_{001})\!\in\!\Bbb{C}^3$ is a line
through the origin.

\subsection{The Degree-Two Number}
\label{n2_subs}

\noindent
The computation of the number $n_2$ is very similar.
A degree-two curve in $\PP$ is described by a nonzero
degree-two homogeneous polynomial
$$s_{a_{2,0,0},a_{1,1,0},a_{1,0,1},
a_{0,2,0},a_{0,1,1},a_{0,0,2}}=\sum_{j+k+l=2}a_{jkl}X^jY^kZ^l.$$
Thus, the space of degree-two curves in $\PP$ is 
$${\cal D}_2=\big\{
\big(a_{2,0,0},a_{1,1,0},a_{1,0,1},a_{0,2,0},a_{0,1,1},a_{0,0,2})
\!\in\!\Bbb{C}^6\!-\!\{0\}\big\}\big/\Bbb{C}^*\approx
\Bbb{P}^5.$$
If $p_i\!=\![X_i,Y_i,Z_i]$ for $i\!=\!1,\ldots,5$ are five points
in $\PP$, the subset of conics that pass through these points
is the set of elements $[(a_{jkl})_{j+k+l=2}]\!\in\!{\cal D}_2$ 
such that
\begin{equation}
\label{n2_e}
\sum_{j+k+l=2}\!\!a_{jkl}X_i^jY_i^kZ_i^l=0
\qquad\hbox{for}\quad i=1,\ldots,5.
\end{equation}
Each of these five linear equations determines a hyperplane $H_i$ 
in~${\cal D}_2$.\\

\noindent
We assume that the five points~$p_i$ do not lie on any pair of lines in~$\PP$.
Then by Lemma~\ref{genus_lmm}, every conic passing through 
the five points~$p_i$ is smooth and of genus~zero.
It follows that any two distinct conics ${\cal C}_1$ and ${\cal C}_2$ 
passing through the five points~$p_i$ must intersect at most
$2\cdot 2\!=\!4$ points; see Lemma~\ref{bezout_lmm}.
Thus, the system~\e_ref{n2_e} of five equations must have at most
one solution~${\cal D}_2$, and such a solution represents 
a plane rational conic through the five points in~$\PP$.
On the other hand, the five hyperplanes~$H_i$ in ${\cal D}_2$
must have at least a point in common, since the poincare dual
of a hyperplane generates~$H^*(\P;\Bbb{Z})$.
In simpler terms, the solution space of the system~\e_ref{n2_e} of
five linear homogeneous equations on~$\Bbb{C}^6$ must contain
a line through the origin. We conclude that $n_2\!=\!1$.

\subsection{The Degree-Three Number}
\label{n3_subs}

\noindent
Computing the number $n_3$ requires a bit more care.
Similarly to the previous two subsections, 
the space of cubics in~$\PP$ is described~by
$${\cal D}_3=
\big\{(a_{jkl})_{j+k+l=3}\!\in\!\Bbb{C}^{10}\!-\!\{0\}\big\}\big/
\Bbb{C}^*\approx\Bbb{P}^9.$$
For a generic $\under{a}\!\in\!{\cal D}_3$,
the section $s_{\under{a}}$ of the bundle 
$\ga^{*\otimes3}\!\lra\!\PP$ is transverse to the zero~set.
Thus, by Lemma~\ref{genus_lmm}, a typical cubic is smooth and 
of genus one, not~zero.\\

\noindent
Let $p_i\!=\![X_i,Y_i,Z_i]$ for $i\!=\!1,\ldots,8$ 
be eight points in $\PP$ that do not lie on the union
of any line and any conic in~$\PP$.
It can then be shown that 
if the cubic ${\cal C}_{\under{a}}$ passes through these eight points,
the section~$s_{\under{a}}$ has at most one singular point.
In such a case, the curve~${\cal C}_{\under{a}}$ is 
a sphere with two points identified.
In other words, a circle on a torus collapses to a~point.
This fact is immediate from the algebraic-geometry point of view,
but can also be checked directly.
Thus, the number~$n_3$ is the number of plane cubics that
pass through the eight points $p_1,\ldots,p_8$ and have a singular point.
This singular point will be a simple node;  see Figure~\ref{sing_fig}
on page~\pageref{sing_fig}.\\

\noindent
As in the previous two subsections, the space $H_i$ of elements 
$\under{a}\!\in\!{\cal D}_3$ such that $p_i\!\in\!{\cal C}_{\under{a}}$
is a hyperplane. With our assumption on the eight points,
the eight hyperplanes intersect transversally, and thus
$${\cal D}\equiv\bigcap_{i=1}^{i=8}H_i
\approx\Bbb{P}^1.$$
In simpler words, the eight equations analogous to~\e_ref{n2_e}
are linearly independent. Thus, the space of solution of
the corresponding system of equations on $\Bbb{C}^{10}$
is a plane through the origin, which corresponds to
a line~$\Bbb{P}^1$ in ${\cal D}_3\!\approx\!\Bbb{P}^9$.\\

\noindent
By the above, we need to determine the cardinality of the set
$${\cal Z}=\big\{\big([\under{a}],x\big)\!\in\!{\cal S}\!:
ds_{\under{a}}\big|_x\!=\!0\big\}, \quad\hbox{where}\quad
{\cal S}=\big\{\big([\under{a}],x\big)\!\in\!{\cal D}\!\times\!\PP\!:
s_{\under{a}}(x)\!=\!0\big\}.$$
An element of the subspace ${\cal S}$ of ${\cal D}\!\times\!\PP$ is
a cubic through the eight points $p_1,\ldots,p_8$ with a choice
of a point on~it.
Such an element $([\under{a}],x)$ lies in~${\cal Z}$ if $s_{\under{a}}$
is not transverse to the zero set at~$x$.\\

\noindent
Let $\pi_0,\pi_1\!:{\cal D}\!\times\!\PP\!\lra\!{\cal D},\PP$
be the two projection maps.
If $\ga_{\cal D}\!\lra\!{\cal D}$ and $\ga_{\PP}\!\lra\!\PP$
are the tautological line bundles, we set 
$$\ga_0\!=\!\pi_0^*\ga_{\cal D}\lra {\cal D}\!\times\!\PP
\qquad\hbox{and}\qquad
\ga_1\!=\!\pi_1^*\ga_{\PP}\lra {\cal D}\!\times\!\PP.$$
A homogeneous polynomial in three variables of degree~$d$ 
induces a section of the bundle \hbox{$\ga^{*\otimes d}\!\!\lra\!\PP$}.
For the same reason, the map
$$\big\{\under{a}\!\in\!\Bbb{C}^2\!:[\under{a}]\!\in\!{\cal D}\big\}
\!\times\!\Bbb{P}^2\lra\ga_{\PP}^{*\otimes3},\qquad
(\under{a},x)\lra s_{\under{a}}(x),$$
induces a section $\psi_0$ of the line bundle 
$\ga_0^*\!\otimes\!\ga_1^{*\otimes3}\!\lra\!{\cal D}\!\times\!\PP$.
This section is transverse to the zero~set.
Thus, ${\cal S}\!=\!\psi_0^{-1}(0)$ is a smooth submanifold of
${\cal D}\!\times\!\PP$; see Lemma~\ref{euler_lmm} below.\\

\noindent
If $([\under{a}],x)\!\in\!{\cal S}$, $s_{\under{a}}(x)\!=\!0$, and
thus $ds_{\under{a}}|_x$ is well-defined.
The~map 
$$\big\{\big([\under{a}],x)\!\in\!\Bbb{C}^2\!\times\!\PP\!:
s_{\under{a}}(x)\!=\!0\big\}\lra
\ga_{\PP}^{*\otimes3}\!\otimes\!T^*\PP,\qquad
(\under{a},x)\lra ds_{\under{a}}\big|_x,$$
induces a section $\psi_1$ of the vector bundle 
$\ga_0^*\!\otimes\!\ga_1^{*\otimes3}\!\otimes\!\pi_1^*T^*\PP\!\lra\!{\cal S}$.
This section is transverse to the zero~set.
Thus, by Lemma~\ref{euler_lmm},
\begin{equation*}\begin{split}
n_3=|{\cal Z}|&=\big|\psi_1^{-1}(0)\big| =
\blr{e\big(\ga_0^*\!\otimes\!\ga_1^{*\otimes3}\!\otimes\pi_1T^*\PP\big),
[{\cal S}]}\\
& =\blr{c_2\big(\ga_0^*\!\otimes\!\ga_1^{*\otimes3}\!\otimes\pi_1^*T^*\PP\big)
\hbox{PD}_{{\cal D}\times\PP}([{\cal S}]),\big[{\cal D}\!\times\!\PP\big]}\\
&=\blr{(3ya\!+\!3a^2)(y\!+\!3a),\big[{\cal D}\!\times\!\PP\big]}=12,
\end{split}\end{equation*}
where $y\!=\!\pi_0^*c_1(\ga_{\cal D}^*)$ and 
$a\!=\!\pi_1^*c_1(\ga_{\PP}^*)$.

\begin{lmm}
\label{euler_lmm}
If $M$ is a compact oriented manifold, 
$V\!\lra\!M$ is an oriented vector bundle, and
$\psi\!\in\!\Ga(M;V)$ is transverse to the zero~set,
the space $\psi^{-1}(0)$ is a smooth oriented submanifold of~$M$ and
$$\hbox{PD}_M\big([\psi^{-1}(0)]\big)=e(V)\in H^*(M;\Bbb{Z}),$$
where $e(V)$ is the euler class of~$V$.
\end{lmm}

\noindent
This lemma is a standard fact in differential topology;
see Sections 9-12 of~\cite{MiSt}.
It implies that if the dimension of~$M$ and the rank of~$V$
are the same, the set $s^{-1}(0)$ is finite and its signed cardinality 
is given~by
$$^{\pm}\!\big|s^{-1}(0)\big|=\big\lan e(V),[M]\big\ran.$$
In fact, this is the only case of Lemma~\ref{euler_lmm} 
we would have needed
if we extended the section~$\psi_1$ over the entire space 
${\cal D}\!\times\!\PP$ by using the canonical connection
of the hermitian holomorphic vector bundle $\ga\!\lra\!\PP$;
see~\cite{GriH}.

\section{The Degree-Four Number}
\label{deg4_sec}

\subsection{Summary}
\label{summary4_subs}

\noindent
In this section we use the general approach of Subsection~\ref{n3_subs}
to compute the number $n_4$.
Since the genus of a smooth plane quartic is three by Lemma~\ref{genus_lmm}, 
we will need to determine the number of quartics that pass through 11 points
in $\PP$ and have three nodes.
This number is one-sixth the cardinality of the~set
$$\tilde{\cal N}_3\equiv \big\{
([\under{a}],x_1,x_2,x_3)\!\in\!{\cal D}\!\times\!\PP_1\times\!\PP_2\!\times\!\PP_3\!:
x_i\!\neq\!x_j~\forall i\!\neq\!j;~
s_{\under{a}}(x_i)=0,~ds_{\under{a}}|_{x_i}\!=\!0
~\forall i\!=\!1,2,3\big\},$$
where ${\cal D}\!\approx\!\Bbb{P}^3$ is the space of quartics 
that pass through the eleven chosen points and $\PP_i\!=\!\PP$.\\

\noindent
Similarly to Subsection~\ref{n3_subs}, each of the sections
$$\vph_i([\under{a}],x_i)=
\big(s_{\under{a}}(x_i),ds_{\under{a}}|_{x_i}\big)
\in \ga_0^*\!\otimes\!\ga_i^{*\otimes4} \oplus 
 \ga_0^*\!\otimes\!\ga_i^{*\otimes4}\!\otimes\!T\PP_i$$
is transverse to the zero set over ${\cal D}\!\times\!\PP_i$.
However, the section
$$\vph\equiv\vph_1\oplus\vph_2\oplus\vph_3$$
is not transverse to the zero set over 
${\cal D}\!\times\!\PP_1\!\times\!\PP_2\!\times\!\PP_3$.
For example, the zero set of $\vph$ contains the two-dimensional space
$$\big\{([\under{a}],x,x,x)\!:s_{\under{a}}(x)\!=\!0,~ ds_{\under{a}}|_x=0\big\}.$$
Thus, $|\tilde{\cal N}_3|$ is not the euler class of the bundle
$$V\equiv \bigoplus_{i=1}^{i=3}\ga_0^*\!\otimes\!\ga_i^{*\otimes4} \oplus 
 \ga_0^*\!\otimes\!\ga_i^{*\otimes4}\!\otimes\!T^*\PP_i
\lra M\equiv{\cal D}\!\times\!\PP_1\!\times\!\PP_2\!\times\!\PP_3.$$
On the other hand, $\vph$ is transverse to the zero set over the ``main stratum"
of~$M$:
$$M^0\equiv
\big\{([\under{a}],x_1,x_2,x_3)\!\in\!M\!: x_i\!\neq\!x_j~\forall i\!\neq\!j\big\}.$$
Thus, $|\tilde{\cal N}_3|$ is the euler class of the bundle $V$ minus
the {\it $\vph$-contribution} to $e(V)$ from the ``boundary" of~$M$:
$$|\tilde{\cal N}_3|= \lr{e(V),M}-{\cal C}_{\partial M}(\vph),
\qquad\hbox{where}\qquad  \partial M=M-M^0.$$
The number ${\cal C}_{\partial M}(\vph)$ is the signed number of zeros
of the bundle section $\vph\!+\!\nu$, for a small generic perturbation~$\nu$,
that lie near $\partial M$.
If $\partial M\!=\!\sqcup_i{\cal Z}_i$ is a stratification of $\partial M$,
$${\cal C}_{\partial M}(\vph)=\sum_i{\cal C}_{{\cal Z}_i}(\vph).$$
If this stratification is sufficiently fine,
each of the numbers ${\cal C}_{{\cal Z}_i}(\vph)$ is a certain multiple
of the number of zeros of an affine bundle map
between vector bundles over~$\bar{\cal Z}_i$.
The latter number can be computed through a reductive procedure,
described in detail in~\cite{Z1} and~\cite{Z2} and implemented
in the relevant cases in Subsections~\ref{p2_subs} and~\ref{p3_subs} below.\\

\begin{table}
\begin{center}
\begin{tabular}{||c|c|c|c||}
\hline\hline
set& singularities& \#pts& card.\\
\hline\hline
${\cal N}_1$& 1 node& 13& 27\\
\hline
${\cal N}_{1,1}$& 1 node on a fixed line& 12& 9\\
\hline
${\cal K}_1$& 1 cusp& 12& 72\\
\hline
${\cal K}_{1,1}$& 1 cusp on a fixed line& 11& 20\\
\hline
${\cal T}_1$& 1 tacnode& 11& 200\\
\hline\hline
${\cal N}_2$& 2 nodes& 12& 225\\
\hline
${\cal N}_{2,1}$& 2 nodes, one on a fixed line& 11& 170\\
\hline
${\cal K}_2$& 1 node and 1 cusp& 11& 840\\
\hline\hline
${\cal N}_3$& 3 nodes& 11& 675\\
\hline\hline
\end{tabular}
\caption{Some Characteristic Numbers of Plane Quartics}
\label{quartics_table}
\end{center}
\end{table}

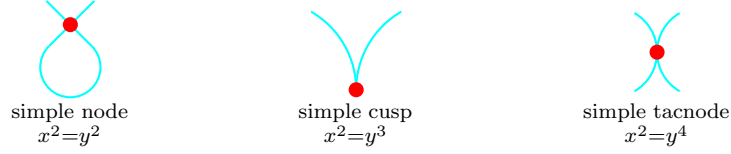
\begin{figure}
\begin{pspicture}(-5,-2.5)(10,.2)
\psset{unit=.2cm}
\psarc[linecolor=cyan](-3,-5){2}{135}{45}
\psline[linecolor=cyan](-4.41,-3.59)(-1.41,-.59)
\psline[linecolor=cyan](-1.59,-3.59)(-4.59,-.59)
\pscircle[fillstyle=solid,fillcolor=red,linecolor=red](-3,-2.17){.5}
\rput(-3,-8){$_{\text{simple node}}$}
\rput(-3,-9.6){$_{x^2=y^2}$}
\psarc[linecolor=cyan](10,-6.5){6}{0}{60}
\psarc[linecolor=cyan](22,-6.5){6}{120}{180}
\pscircle[fillstyle=solid,fillcolor=red,linecolor=red](16,-6.5){.5}
\rput(16,-8){$_{\text{simple cusp}}$}
\rput(16,-9.6){$_{x^2=y^3}$}
\psarc[linecolor=cyan](33,-4){3}{-60}{60}
\psarc[linecolor=cyan](39,-4){3}{120}{240}
\pscircle[fillstyle=solid,fillcolor=red,linecolor=red](36,-4){.5}
\rput(36,-8){$_{\text{simple tacnode}}$}
\rput(36,-9.6){$_{x^2=y^4}$}
\end{pspicture}
\caption{Simple Node, Simple Cusp, and Simple Tacnode}
\label{sing_fig}
\end{figure}

\noindent
In order to simplify the computation of $|\tilde{\cal N}_3|$,
we will essentially be adding one point at a time.
This computation will require knowing the numbers
of plane quartics with various one- and two-point singularities.
These numbers, along with $|{\cal N}_3|$, are given in Table~\ref{quartics_table}.
For example, according to this table,
the cardinality of the set ${\cal N}_{2,1}$ of plane quartics
that pass through 11 points in general position and 
have two nodes, one of which lies on a fixed general line, is~$170$.
Figure~\ref{sing_fig} shows a simple node, a simple cusp, and a simple tacnode.
If $s$ is a section of $\ga^{*\otimes d}$ and $x\!\in\!s^{-1}(0)$
is a node of $s^{-1}(0)$, then $ds|_x\!=\!0$.
We describe the analogous cuspidal and tacnodal condition on $s$ 
in the next subsection.
All numbers in Table~\ref{quartics_table} are computed 
in Subsections~\ref{p1_subs}-\ref{p3_subs}.\\

\noindent
Finally, we note that a plane quartic that has 3 nodes and passes 
through 11 points is either irreducible, in which case it is rational,
or a union of a smooth cubic, passing through 9 of the points,
and a line, passing through the remaining 2 points.
By the same argument as in Subsections~\ref{n1_subs} and~\ref{n2_subs},
the number of plane cubics passing through 9 points in general position is~1.
Thus, by the last row of Table~\ref{quartics_table}, the number of rational
quartics passing through 11 points in general position in $\PP$~is
$$n_4=675-\binom{11}{2}\cdot 1\cdot 1=620.$$\\

\noindent
The computations of Subsections~\ref{p1_subs}-\ref{p3_subs} generalize easily to 
plane curves of arbitrary degree, essentially by replacing $4$ by~$d$ everywhere.
The results and the arguments are summarized in Table~\ref{curvegen_table} 
and in Subsection~\ref{curvgen_subs}, respectively.

\subsection{Quartics with One Singular Point}
\label{p1_subs}

\noindent
Throughout the rest of Section~\ref{deg4_sec},
we denote by $p_1,\ldots,p_{13}$ thirteen points in general position in~$\PP$
and by ${\cal D}_4\!\approx\!\Bbb{P}^{14}$ the space of plane quartics.
In this subsection, we compute the first five numbers in Table~\ref{quartics_table}.

\begin{lmm}
\label{node1_lmm}
The number $|{\cal N}_1|$ of plane quartics that have a node and pass
through $13$ points in general position is~$27$.
The number $|{\cal N}_{1,1}|$ of plane quartics that have a node 
on a fixed general line and pass through $12$ points in general position is~$9$.
\end{lmm}

\noindent
{\it Proof:}
(1) Let ${\cal D}\!\approx\!\Bbb{P}^1\!\subset\!{\cal D}_4$ denote 
the subspace of plane quartics that pass through the points $p_1,\ldots,p_{13}$.
With notation as in Subsection~\ref{n3_subs}, let
\begin{gather*}
{\cal N}_1=\big\{([\under{a}],x)\!\in\!{\cal D}\!\times\!\Bbb{P}^2\!:
\vph([\under{a}],x)=0\big\},
\qquad\hbox{where}\\
\vph\in\Ga\big({\cal D}\!\times\!\PP;
\ga_0^*\!\otimes\!\ga_1^{*\otimes4}\oplus
\ga_0^*\!\otimes\!\ga_1^{*\otimes4}\!\otimes\!T^*\PP\big),\qquad
\vph([\under{a}],x)=\big(s_{\under{a}}(x),ds_{\under{a}}|_x\big).
\end{gather*}
Since the section $\vph$ is transverse to the zero set, by Lemma~\ref{euler_lmm},
\begin{equation*}\begin{split}
|{\cal N}_1|=|\vph^{-1}(0)|&
=\blr{e\big(\ga_0^*\!\otimes\!\ga_1^{*\otimes4}\oplus
\ga_0^*\!\otimes\!\ga_1^{*\otimes4}\!\otimes\!T^*\PP\big),
{\cal D}\!\times\!\PP}\\
& =\blr{c_1\big(\ga_0^*\!\otimes\!\ga_1^{*\otimes4}\big)
c_2\big(\ga_0^*\!\otimes\!\ga_1^{*\otimes4}\!\otimes\!T^*\PP\big),
{\cal D}\!\times\!\PP}\\
&=\blr{(y\!+\!4a)(y^2\!+\!5ya\!+\!7a^2),{\cal D}\!\times\!\PP}=27.
\end{split}\end{equation*}
(2) Let ${\cal D}\!\approx\!\PP\!\subset\!{\cal D}_4$ denote 
the subspace of plane quartics that pass through the points $p_1,\ldots,p_{12}$.
Let $\Bbb{P}^1\!\subset\!\PP$ be a general line in~$\PP$.
We~put
\begin{gather*}
{\cal N}_{1,1}=\big\{([\under{a}],x)\!\in\!{\cal D}\!\times\!\Bbb{P}^1\!:
\vph([\under{a}],x)=0\big\},
\qquad\hbox{where}\\
\vph\in\Ga\big({\cal D}\!\times\!\Bbb{P}^1;
\ga_0^*\!\otimes\!\ga_1^{*\otimes4}\oplus
\ga_0^*\!\otimes\!\ga_1^{*\otimes4}\!\otimes\!T^*\PP|_{\Bbb{P}^1}\big),\qquad
\vph([\under{a}],x)=\big(s_{\under{a}}(x),ds_{\under{a}}|_x\big).
\end{gather*}
Since the section $\vph$ is transverse to the zero set, by Lemma~\ref{euler_lmm},
\begin{equation*}\begin{split}
|{\cal N}_{1,1}|=|\vph^{-1}(0)|&
=\blr{e\big(\ga_0^*\!\otimes\!\ga_1^{*\otimes4}\oplus
\ga_0^*\!\otimes\!\ga_1^{*\otimes4}\!\otimes\!T^*\PP\big),
{\cal D}\!\times\!\Bbb{P}^1}\\
& =\blr{c_1\big(\ga_0^*\!\otimes\!\ga_1^{*\otimes4}\big)
c_2\big(\ga_0^*\!\otimes\!\ga_1^{*\otimes4}\!\otimes\!T^*\PP\big),
{\cal D}\!\times\!\Bbb{P}^1}\\
&=\blr{(y\!+\!4a)(y^2\!+\!5ya\!+\!7a^2),{\cal D}\!\times\!\Bbb{P}^1}=9.
\end{split}\end{equation*}

\begin{lmm}
\label{cusp1_lmm}
The number $|{\cal K}_1|$ of plane quartics that have a cusp and pass
through $12$ points in general position is~$72$.
The number $|{\cal K}_{1,1}|$ of plane quartics that have a cusp 
on a fixed general line and pass through $11$ points in general position is~$20$.
\end{lmm}

\noindent
{\it Proof:}
(1) Let ${\cal D}\!\approx\!\Bbb{P}^2$ be as in (2) of the proof of 
Lemma~\ref{node1_lmm}. We put
$${\cal N}_1'=\big\{([\under{a}],x)\!\in\!{\cal D}\!\times\!\PP\!:
s_{\under{a}}(x)\!=\!0,~ds_{\under{a}}|_x\!=\!0\big\}.$$
If $([\under{a}],x)\!\in\!{\cal N}_1'$, we denote by
$$H_{\under{a},x}\in
\Ga\big({\cal N}_1';
\Hom(T\PP,\ga_0^*\!\otimes\!\ga_1^{*\otimes4}\!\otimes\!T^*\PP)\big)$$
the Hessian of $s_{\under{a}}$ at $x$, 
i.e.~the total second derivative of $s_{\under{a}}$ at~$x$.
Let
\begin{gather*}
{\cal K}_1=\big\{([\under{a}],x)\!\in\!{\cal N}_1'\!:\vph([\under{a}],x)\!=\!0\big\}
\qquad\hbox{where}\\
\vph\in\Ga\big({\cal N}_1';
(\ga_0^*\!\otimes\!\ga_1^{*\otimes4}\!\otimes\!\La^2T^*\PP)^{\otimes2}\big),
\quad \vph([\under{a}],x)=\hbox{det}\,H_{\under{a},x}.
\end{gather*}
Since the section $\vph$ is transverse to the zero set, by Lemma~\ref{euler_lmm},
\begin{equation*}\begin{split}
|{\cal K}_1|=|\vph^{-1}(0)|&
=\blr{e\big((\ga_0^*\!\otimes\!\ga_1^{*\otimes4}\!\otimes\!\La^2T^*\PP)^{\otimes2}\big),
{\cal N}_1'}\\
&=2\blr{y\!+\!a,{\cal N}_1'}=2\, \big(|{\cal N}_1|\!+\!|{\cal N}_{1,1}|\big)
=2\, (27\!+\!9) =72.
\end{split}\end{equation*}
(2) Similarly, let ${\cal D}\!\approx\!\Bbb{P}^3\!\subset\!{\cal D}_4$ denote 
the subspace of plane quartics that pass through the points $p_1,\ldots,p_{11}$.
Let $\Bbb{P}^1\!\subset\!\PP$ be a general line in~$\PP$.
We~put
\begin{alignat*}{1}
{\cal N}_{1,1}'&=\big\{([\under{a}],x)\!\in\!{\cal D}\!\times\!\Bbb{P}^1\!:
s_{\under{a}}(x)\!=\!0,~ds_{\under{a}}|_x\!=\!0\big\};\\
{\cal K}_{1,1}&=\big\{([\under{a}],x)\!\in\!{\cal N}_{1,1}'\!:
\hbox{det}\,H_{\under{a},x}=0\big\}.
\end{alignat*}
Then, by Lemma~\ref{euler_lmm},
\begin{equation*}\begin{split}
|{\cal K}_{1,1}|&
=\blr{e\big((\ga_0^*\!\otimes\!\ga_1^{*\otimes4}\!\otimes\!\La^2T^*\PP)^{\otimes2}\big),
{\cal N}_{1,1}'}\\
&=2\blr{y\!+\!a,{\cal N}_{1,1}'}
=2\, \big(|{\cal N}_{1,1}|\!+\!\lr{a,{\cal N}_{1,1}'}\big)
=2\, (9\!+\!1) =20.
\end{split}\end{equation*}
Note the number $\lr{a,{\cal N}_{1,1}'}$ of plane quartics that pass
through 11 points and have a node at a fixed twelfth point is~$1$,
since all conditions on $\under{a}\!\in\!{\cal D}_4$ are linear,
as in Subsections~\ref{n1_subs} and~\ref{n2_subs}.

\begin{lmm}
\label{tacnode1_lmm}
The number $|{\cal T}_1|$ of plane quartics that have a tacnode and pass
through $11$ points in general position is~$200$.
\end{lmm}

\noindent
{\it Proof:}
Let ${\cal D}\!\approx\!\Bbb{P}^3$ be as in (2) of the proof of 
Lemma~\ref{cusp1_lmm}. We put
$${\cal N}_1''=\big\{([\under{a}],x)\!\in\!{\cal D}\!\times\!\PP\!:
s_{\under{a}}(x)\!=\!0,~ds_{\under{a}}|_x\!=\!0\big\}, 
\qquad M=\Bbb{P}T\PP|_{{\cal N}_1''}.$$
We denote by $\ga\!\lra\!M$ the tautological line bundle and by
$$\tilde{H}_{\cdot,\cdot}\in 
\Ga\big(M;\Hom(\ga,\ga_0^*\!\otimes\!\ga_1^{*\otimes4}\!\otimes\!T^*\PP)\big)$$
the bundle map induced by $H_{\cdot,\cdot}$.
Let
\begin{gather*}
{\cal K}_1'=\big\{([\under{a}],x)\!\in\!M: \tilde{H}_{\under{a},x}=0\big\},
\qquad
{\cal T}_1=\big\{([\under{a}],x)\!\in\!{\cal K}_1'\!:
\vph(\under{a},x)=0\big\},\\
\hbox{where}\qquad
\vph\in\Ga\big(M;\Hom(\ga^{\otimes3},\ga_0^*\!\otimes\!\ga_1^{*\otimes4})\big),
\quad
\vph([\under{a}],x)={\cal D}_{\under{a},x}^3,
\end{gather*}
and ${\cal D}_{\under{a},x}^3$ is the third derivative of $s_{\under{a}}$ at $x$.
Let $\la\!=\!c_1(\ga^*)$.
Since the sections $\vph$ and $\tilde{H}_{\cdot,\cdot}$ are transverse 
to the zero set, 
by Lemma~\ref{euler_lmm},
\begin{equation*}\begin{split}
|{\cal T}_1|=|\vph^{-1}(0)|
&=\blr{e\big(\ga^{*\otimes3}\!\otimes\!\ga_0^*\!\otimes\!\ga_1^{*\otimes4}\big),{\cal K}_1'}\\
&=\blr{e\big(\ga^{*\otimes3}\!\otimes\!\ga_0^*\!\otimes\!\ga_1^{*\otimes4}\big)
e\big(\ga^*\!\otimes\!\ga_0^*\!\otimes\!\ga_1^{*\otimes4}\!\otimes\!T^*\PP\big),M}\\
&=\blr{3\la^3+(7y\!+\!19a)\la^2+(5y^2\!+\!28ya\!+\!41a^2)\la,M}\\
&=\blr{5y^2+7ya+2a^2,{\cal N}_1''}\\
&=5|{\cal N}_1|+7|{\cal N}_{1,1}|+2\lr{a,{\cal N}_{1,1}'}
=5\cdot 27+7\cdot 9+2\cdot 1=200.
\end{split}\end{equation*}

\subsection{Quartics with Two Singular Points}
\label{p2_subs}

\noindent
In this subsection, we compute the three numbers of Table~\ref{quartics_table}
that involve two-point singularities.
As the relevant bundle sections are no longer transverse everywhere,
each of these numbers is the euler class of the corresponding vector bundle
minus the contribution from the "boundary" for the given bundle section.\\

\noindent
Suppose $E,V\!\lra\!M$ are vector bundle such that 
$\hbox{dim}\,M\!+\!\hbox{rk}\,E\!=\!\hbox{rk}\,V$ and
$$\al\in\Ga\big(M;\Hom(E,V)\big).$$
If $\nu\!\in\!\Ga(M;V)$ is a generic section, the affine bundle map
$$\psi_{\al,\nu}\!:E\lra V, \qquad  \psi_{\al,\nu}(m;e)=\al(m;e)+\nu(m),$$
has a finite number of transverse zeros.
By Lemma~3.14 in~\cite{Z1} and Proposition~2.18A in~\cite{Z2}, the signed cardinality 
of $\psi_{\al,\nu}^{-1}(0)$ is independent of the choice of $\nu$.
We denote this cardinality by~$N(\al)$.

\begin{lmm}
\label{node2_lmm1}
The number $|{\cal N}_2|$ of plane quartics that have two nodes and pass
through $12$ points in general position is~$225$.
The number $|{\cal N}_{2,1}|$ of plane quartics that have two nodes,
one of which lies on a fixed general line,
 and pass through $11$ points in general position is~$170$.
\end{lmm}

\noindent
{\it Proof:}
(1) Let ${\cal N}_1'\!\subset\!{\cal D}\!\times\!\PP_1$ be defined
as in (1) of the proof of Lemma~\ref{cusp1_lmm}.
We~put
\begin{gather*}
M={\cal N}_1'\!\times\!\PP_2, \quad
M^0=\big\{([\under{a}],x_1,x_2)\!\in\!M\!:x_1\!\neq\!x_2\},\quad
\partial M=M-M^0, \quad
\tilde{\cal N}_2=\vph^{-1}(0)\cap M^0,\\
\hbox{where}\qquad
\vph\!\in\!\Ga(M;\ga_0^*\!\otimes\!\ga_2^{*\otimes4}\oplus
\ga_0^*\!\otimes\!\ga_2^{*\otimes4}\!\otimes\!T^*\PP_2),  \quad
\vph([\under{a}],x_1,x_2)=\big(s_{\under{a}}(x_2),ds_{\under{a}}|_{x_2}\big),
\quad \ga_2\!=\!\pi_2^*\ga_{\PP_2},
\end{gather*}
and $\pi_2\!:M\!\lra\!\PP_2$ is the projection onto the last component.
Since $\vph|_{M^0}$ is transverse to the zero set,
\begin{equation}\label{node2_e1}\begin{split}
|\tilde{\cal N}_2|=\, ^{\pm}\big|\vph^{-1}(0)\cap M^0\big|
&= \blr{e(\ga_0^*\!\otimes\!\ga_2^{*\otimes4}\oplus
\ga_0^*\!\otimes\!\ga_2^{*\otimes4}\!\otimes\!T^*\PP_2),M}
-{\cal C}_{\partial M}(\vph)\\
&= \blr{(y\!+\!4a_2)(y^2\!+\!5ya_2\!+\!7a_2^2),{\cal N}_1'\!\times\!\PP_2}
-{\cal C}_{\partial M}(\vph)\\
&= 27\lr{y,{\cal N}_1'}-{\cal C}_{\partial M}(\vph)
=27|{\cal N}_1|-{\cal C}_{\partial M}(\vph)=27\cdot 27-{\cal C}_{\partial M}(\vph),
\end{split}\end{equation}
where $a_2\!=\!\pi_2^*c_1(\ga_{\PP_2}^*)$.
In order to determine ${\cal C}_{\partial M}(\vph)$,
we split $\partial M$ into two strata:
$${\cal Z}_1=\big\{([\under{a}],x,x)\!: ([\under{a}],x)\!\in\!{\cal N}_1'\!-\!{\cal K}_1\big\},
\qquad
{\cal Z}_0=\big\{([\under{a}],x,x)\!: ([\under{a}],x)\!\in\!{\cal K}_1\big\}.$$
With appropriate identifications, for some $C\!\in\!C({\cal N}_1';\R^+)$, 
\begin{equation}\label{node2_e2}
\big|\vph([\under{a}],x,v)-H_{\under{a},x}v\big|\le C([\under{a}],x)|v|^2
\quad\forall\, ([\under{a}],x,x)\!\in\!\partial M, ~
v\!\in\!\hbox{Norm}_M\partial M\big|_{([\under{a}],x,x)} \approx T_x\PP_1.
\end{equation}
By definition of the set ${\cal K}_1$,
\begin{equation}\label{node2_e3}
\big|H_{\under{a},x}v\big|\ge C([\under{a}],x)^{-1}|v|
\quad\forall\, ([\under{a}],x)\!\in\!{\cal N}_1'\!-\!{\cal K}_1, \, v\!\in\! T_x\PP_1.
\end{equation}
By \e_ref{node2_e2}, \e_ref{node2_e3}, and a rescaling and cobordism argument
as in Subsection~3.1 of~\cite{Z1},
\begin{gather}\label{node2_e4}
{\cal C}_{{\cal Z}_1}(\vph)=N(\al),   \qquad\hbox{where}\\
\al\in\Ga\big({\cal N}_1';\Hom(T\PP,\ga_0^*\!\otimes\!\ga_1^{*\otimes4}\oplus
\ga_0^*\!\otimes\!\ga_1^{*\otimes4}\!\otimes\!T^*\PP)\big),
\quad \al([\under{a}],x;v)=(0,H_{\under{a},x}v).  \notag
\end{gather}
On the other hand, suppose $([\under{a}],x)\!\in\!{\cal K}_1$.
We denote by ${\cal L}_{(\under{a},x)}\!\subset\!T\PP$ 
the kernel of $H_{\under{a},x}$ and by  ${\cal L}_{(\under{a},x)}^{\perp}$
its orthogonal complement. 
Let $N_{(\under{a},x)}$ be the normal bundle of ${\cal K}_1$ in
${\cal N}_1'$ at $([\under{a}],x)$.
Then, with appropriate identifications, for some $\be_2,\be_3\!\in\!\Bbb{C}^*$,
$\be_4\!\in\!\Bbb{C}$, and $C\!\in\!\Bbb{R}^+$,
\begin{gather}\label{node2_e5}
\big|\vph([\under{a}],x;u,v,w)-\al_0(u,v,w)\big| \le C\big(|v|^4\!+\!|w|^2)
\quad\forall\, u\!\in\!N_{(\under{a},x)},\,
v\!\in\!{\cal L}_{(\under{a},x)},\, w\!\in\!{\cal L}_{(\under{a},x)}^{\perp},\\
\hbox{where}\quad
\al_0(u,v,w)=\big(\frac{1}{2}uv^2+\frac{1}{3}\be_3v^3,uv+\be_3v^2+\be_4v^3,\be_2w).
\notag
\end{gather}
Here $\be_2$ is the second derivative of $s_{\under{a}}$ at $x$ 
along ${\cal L}_{(\under{a},x)}^{\perp}$ and
$2\be_3$ is the third derivative of $s_{\under{a}}$ at $x$ 
along~${\cal L}_{(\under{a},x)}$.
Since the polynomial~$\al_0$ is three-to-one near the origin,
it follows from \e_ref{node2_e5} that each point of ${\cal Z}_0\!\approx\!{\cal K}_1$
contributes~$3$ to ${\cal C}_{{\cal Z}_0}(\vph)$.
From~\e_ref{node2_e4} and Lemmas~\ref{cusp1_lmm} and~\ref{node2_lmm2}, we conclude
that
\begin{equation}\label{node2_e6}
{\cal C}_{\partial M}(\vph)
={\cal C}_{{\cal Z}_1}(\vph)+{\cal C}_{{\cal Z}_0}(\vph)
=63+3|{\cal K}_1|=63+ 3\cdot 72=279.
\end{equation}
The first claim of the lemma follows from \e_ref{node2_e1} and~\e_ref{node2_e6},
since ${\cal N}_2\!=\!\tilde{\cal N}_2/S_2$, where $S_2$ is the symmetric group
on two elements.\\
(2) Similarly, let ${\cal N}_{1,1}'\!\subset\!{\cal D}\!\times\!\PP_1$ be defined
as in (2) of the proof of Lemma~\ref{cusp1_lmm}.
We~put
\begin{gather*}
M={\cal N}_{1,1}'\!\times\!\PP_2, ~~
M^0=\big\{([\under{a}],x_1,x_2)\!\in\!M\!:x_1\!\neq\!x_2\}, ~~
\partial M=M\!-\!M^0, ~~
{\cal N}_{2,1}=\vph^{-1}(0)\cap M^0,\\
\hbox{where}\qquad
\vph\!\in\!\Ga(M;\ga_0^*\!\otimes\!\ga_2^{*\otimes4}\oplus
\ga_0^*\!\otimes\!\ga_2^{*\otimes4}\!\otimes\!T^*\PP_2),  \quad
\vph([\under{a}],x_1,x_2)=\big(s_{\under{a}}(x_2),ds_{\under{a}}|_{x_2}\big).
\end{gather*}
Since $\vph|_{M^0}$ is transverse to the zero set,
\begin{equation}\label{node2_e11}\begin{split}
|{\cal N}_{2,1}|=\, ^{\pm}\big|\vph^{-1}(0)\cap M^0\big|
&= \blr{e(\ga_0^*\!\otimes\!\ga_2^{*\otimes4}\oplus
\ga_0^*\!\otimes\!\ga_2^{*\otimes4}\!\otimes\!T^*\PP_2),M}
-{\cal C}_{\partial M}(\vph)\\
&= 27\lr{y,{\cal N}_{1,1}'}-{\cal C}_{\partial M}(\vph)
=27|{\cal N}_{1,1}|-{\cal C}_{\partial M}(\vph)=27\cdot 9-{\cal C}_{\partial M}(\vph).
\end{split}\end{equation}
We split $\partial M$ into two strata:
$${\cal Z}_1=\big\{([\under{a}],x,x)\!:  
([\under{a}],x)\!\in\!{\cal N}_{1,1}'\!-\!{\cal K}_{1,1}\big\},
\qquad
{\cal Z}_0=\big\{([\under{a}],x,x)\!: ([\under{a}],x)\!\in\!{\cal K}_{1,1}\big\}.$$
By the same argument as in (1) above,
\begin{gather*}
{\cal C}_{{\cal Z}_1}(\vph)=N(\al),   \qquad\hbox{where}\\
\al\in\Ga\big({\cal N}_{1,1}';\Hom(T\PP,\ga_0^*\!\otimes\!\ga_1^{*\otimes4}\oplus
\ga_0^*\!\otimes\!\ga_1^{*\otimes4}\!\otimes\!T^*\PP)\big),
\quad \al([\under{a}],x;v)=(0,H_{\under{a},x}v),
\end{gather*}
while ${\cal C}_{{\cal Z}_0}(\vph)\!=\!3|{\cal K}_{1,1}|$.
Using Lemmas~\ref{cusp1_lmm} and~\ref{node2_lmm2}, we conclude that
\begin{equation}\label{node2_e16}
{\cal C}_{\partial M}(\vph)
={\cal C}_{{\cal Z}_1}(\vph)+{\cal C}_{{\cal Z}_0}(\vph)
=13+ 3\cdot 20=73.
\end{equation}
The second claim of the lemma follows immediately from \e_ref{node2_e11} and~\e_ref{node2_e16}.

\begin{lmm}
\label{node2_lmm2}
If ${\cal N}_1'\!\subset\!{\cal D}\!\times\!\PP$ is
as in (1) of the proof of Lemma~\ref{cusp1_lmm} and 
$$\al\in\Ga\big({\cal N}_1';\Hom(T\PP,\ga_0^*\!\otimes\!\ga_1^{*\otimes4}\oplus
\ga_0^*\!\otimes\!\ga_1^{*\otimes4}\!\otimes\!T^*\PP)\big),
\quad \al([\under{a}],x;v)=(0,H_{\under{a},x}v),$$ 
then $N(\al)\!=\!63$.
If ${\cal N}_{1,1}'\!\subset\!{\cal D}\!\times\!\PP$ is
as in (2) of the proof of Lemma~\ref{cusp1_lmm} and 
$$\al\in\Ga\big({\cal N}_{1,1}';\Hom(T\PP,\ga_0^*\!\otimes\!\ga_1^{*\otimes4}\oplus
\ga_0^*\!\otimes\!\ga_1^{*\otimes4}\!\otimes\!T^*\PP)\big),
\quad \al([\under{a}],x;v)=(0,H_{\under{a},x}v),$$ 
then $N(\al)\!=\!13$.
\end{lmm}

\noindent
{\it Proof:}
(1) We put 
$$M=\Bbb{P}T\PP|_{{\cal N}_1'}, \qquad
\partial M=\big\{([\under{a}],x)\!\in\!M\!: \tilde{H}_{\under{a},x}\!=\!0\big\}
\approx{\cal K}_1,$$
where $\tilde{H}_{\cdot,\cdot}$ is as in the proof of Lemma~\ref{tacnode1_lmm}.
Let
$$\tilde{\al}=(0,\tilde{H})
\in \Ga\big(M;\Hom(\ga,\ga_0^*\!\otimes\!\ga_1^{*\otimes4}\oplus
\ga_0^*\!\otimes\!\ga_1^{*\otimes4}\!\otimes\!T^*\PP)\big)$$
be the section induced by $\al$.
By Lemma~3.14 in~\cite{Z1} or Proposition~2.18A in~\cite{Z2},
\begin{equation}\label{node2_e21}\begin{split}
N(\al)&=\blr{c\big(\ga_0^*\!\otimes\!\ga_1^{*\otimes4}\oplus
\ga_0^*\!\otimes\!\ga_1^{*\otimes4}\!\otimes\!T^*\PP\big)c(T\PP)^{-1},{\cal N}_1'}
-{\cal C}_{\tilde{\al}^{-1}(0)}(\tilde{\al}^{\perp})\\
&=\lr{3y+6a,{\cal N}_1'}-{\cal C}_{\partial M}(\tilde{\al}^{\perp})
=\big(3|{\cal N}_1|+6|{\cal N}_{1,1}|)-{\cal C}_{\partial M}(\tilde{\al}^{\perp}),
\end{split}\end{equation}
where $\tilde{\al}^{\perp}$ is the composition of the linear bundle map $\tilde{\al}$
with the quotient projection map
$$\ga_0^*\!\otimes\!\ga_1^{*\otimes4}\oplus
\ga_0^*\!\otimes\!\ga_1^{*\otimes4}\!\otimes\!T^*\PP
\lra
\big(\ga_0^*\!\otimes\!\ga_1^{*\otimes4}\oplus
\ga_0^*\!\otimes\!\ga_1^{*\otimes4}\!\otimes\!T^*\PP\big)/\Bbb{C}\nu,$$
for a generic nonvanishing section $\nu$.
The claim~\e_ref{node2_e21} can in fact be easily seen directly from the definition
of $N(\al)$.
Since the section $\tilde{H}$ is transverse to the zero set, so is
the section~$\tilde{\al}^{\perp}$ if $\nu$ is generic.
Thus,
\begin{equation}\label{node2_e22}
{\cal C}_{\tilde{\al}^{-1}(0)}(\tilde{\al}^{\perp})= \,
^{\pm}\big|\tilde{\al}^{-1}(0)\big|=|{\cal K}_1|.
\end{equation}
The first claim of the lemma follows from \e_ref{node2_e21} and \e_ref{node2_e22},\
along with Lemmas~\ref{node1_lmm} and~\ref{cusp1_lmm}.\\
(2) Similarly, we put 
\begin{gather*}
M=\Bbb{P}T\PP|_{{\cal N}_{1,1}'}, \qquad
\partial M=\big\{([\under{a}],x)\!\in\!M\!: \tilde{H}_{\under{a},x}\!=\!0\big\}
\approx{\cal K}_{1,1},\\
\tilde{\al}=(0,\tilde{H})
\in \Ga\big(M;\Hom(\ga,\ga_0^*\!\otimes\!\ga_1^{*\otimes4}\oplus
\ga_0^*\!\otimes\!\ga_1^{*\otimes4}\!\otimes\!T^*\PP)\big).
\end{gather*}
By Lemma~3.14 in~\cite{Z1} or Proposition~2.18A in~\cite{Z2},
\begin{equation}\label{node2_e26}\begin{split}
N(\al)&=\blr{c\big(\ga_0^*\!\otimes\!\ga_1^{*\otimes4}\oplus
\ga_0^*\!\otimes\!\ga_1^{*\otimes4}\!\otimes\!T^*\PP\big)c(T\PP)^{-1},{\cal N}_{1,1}'}
-{\cal C}_{\tilde{\al}^{-1}(0)}(\tilde{\al}^{\perp})\\
&=\lr{3y+6a,{\cal N}_{1,1}'}-{\cal C}_{\partial M}(\tilde{\al}^{\perp})
=\big(3|{\cal N}_{1,1}|+6\lr{a,{\cal N}_{1,1}})-{\cal C}_{\partial M}(\tilde{\al}^{\perp}).
\end{split}\end{equation}
As in (1), $\tilde{\al}^{\perp}$ is transverse to the zero, and thus
\begin{equation}\label{node2_e27}
{\cal C}_{\tilde{\al}^{-1}(0)}(\tilde{\al}^{\perp})= \,
^{\pm}\big|\tilde{\al}^{-1}(0)\big|=|{\cal K}_{1,1}|.
\end{equation}
The second claim of the lemma follows from \e_ref{node2_e26} and \e_ref{node2_e27},\
along with Lemmas~\ref{node1_lmm} and~\ref{cusp1_lmm}.

\begin{lmm}
\label{cusp2_lmm1}
The number $|{\cal K}_2|$ of plane quartics that have one node and 
one cusp and pass through $11$ points in general position is~$840$.
\end{lmm}

\noindent
{\it Proof:} Let ${\cal N}_1''\!\subset\!{\cal D}\!\times\!\PP_1$ and
${\cal K}_1'\!\subset\!\Bbb{P}T\PP_1|_{{\cal N}_1''}$ be as in the proof
of Lemma~\ref{tacnode1_lmm}.
We denote by 
$$\tilde{\pi}_1\!:\Bbb{P}T\PP_1|_{{\cal N}_1''}\lra\PP_1$$
the composition of the bundle projection $\Bbb{P}T\PP_1|_{{\cal N}_1''}\!\lra\!{\cal N}_1''$
with $\pi_1$.
We~put
\begin{gather*}
M\!=\!{\cal K}_1'\!\times\!\PP_2, ~~
M^0=\big\{([\under{a}],x_1,x_2)\!\in\!M\!:\tilde{\pi}_1([\under{a}],x_1)\!\neq\!x_2\}, 
~~ \partial M\!=\!M\!-\!M^0, ~~
{\cal K}_2\!=\!\vph^{-1}(0)\cap M^0,\\
\hbox{where}\quad
\vph\!\in\!\Ga(M;\ga_0^*\!\otimes\!\ga_2^{*\otimes4}\oplus
\ga_0^*\!\otimes\!\ga_2^{*\otimes4}\!\otimes\!T^*\PP_2), ~~~
\vph([\under{a}],x_1,x_2)=\big(s_{\under{a}}(x_2),ds_{\under{a}}|_{x_2}\big),
~~~ \ga_2\!=\!\pi_2^*\ga_{\PP_2}.
\end{gather*}
Since $\vph|_{M^0}$ is transverse to the zero set, similarly to~\e_ref{node2_e1},
\begin{equation}\label{cusp2_e1}\begin{split}
|{\cal K}_2|=\, ^{\pm}\big|\vph^{-1}(0)\cap M^0\big|
&= \blr{(y\!+\!4a_2)(y^2\!+\!5ya_2\!+\!7a_2^2),{\cal K}_1'\!\times\!\PP_2}
-{\cal C}_{\partial M}(\vph)\\
&= 27\lr{y,{\cal K}_1'}-{\cal C}_{\partial M}(\vph)
=27|{\cal K}_1|-{\cal C}_{\partial M}(\vph)=27\cdot 72-{\cal C}_{\partial M}(\vph).
\end{split}\end{equation}
We split $\partial M$ into two strata:
$${\cal Z}_1=\big\{([\under{a}],x,\tilde{\pi}_1([\under{a}],x))\!:  
([\under{a}],x)\!\in\!{\cal K}_1'\!-\!{\cal T}_1\big\},
\qquad
{\cal Z}_0=\big\{([\under{a}],x,\tilde{\pi}_1([\under{a}],x))\!:  
([\under{a}],x)\!\in\!{\cal T}_1\big\}.$$
Let $\ga^{\perp}\!\lra\!{\cal K}_1'$ be the orthogonal complement of $\ga$ in $\pi^*T\PP$.
We define the bundle~map
\begin{gather*}
\al\in\Ga\big({\cal K}_1';
\Hom(\ga^{\otimes2}\!\oplus\!\ga^{\perp},\ga_0^*\!\otimes\!\ga_1^{*\otimes4}\oplus
\ga_0^*\!\otimes\!\ga_1^{*\otimes4}\!\otimes\!T^*\PP_1)\big)
\qquad\hbox{by}\\
\al(\tilde{v},w)=\big(0,\frac{1}{2}D_{\under{a},x}^3\tilde{v},\tilde{H}_{\under{a},x}w\big)
\in \ga_0^*\!\otimes\!\ga_1^{*\otimes4}\oplus
\ga_0^*\!\otimes\!\ga_1^{*\otimes4}\!\otimes\!\ga^*
\oplus\ga_0^*\!\otimes\!\ga_1^{*\otimes4}\!\otimes\!\ga^{\perp*}.
\end{gather*}
Note that by definition of the set ${\cal T}_1$,
for some $C\!\in\!C({\cal K}_1';\R^+)$, 
\begin{equation}\label{cusp2_e2}
\big|\al_{[\under{a}],x}(\tilde{v},w)\big|\ge 
C([\under{a}],x)^{-1}\big(|\tilde{v}|\!+\!|w|\big)
\quad\forall\, ([\under{a}],x)\!\in\!{\cal K}_1'\!-\!{\cal T}_1, \, 
(\tilde{v},w)\!\in\!\big(\ga^{\otimes2}\!\oplus\!\ga^{\perp})|_{([\under{a}],x)}.
\end{equation}
On the other hand, with appropriate identifications, 
\begin{gather}\label{cusp2_e3}
\big|\vph([\under{a}],x,v,w)-\al_{[\under{a}],x}(v^{\otimes2},w)\big|
\le  C([\under{a}],x)(|v|^3\!+\!|w|^2\big)\\
\forall~ ([\under{a}],x,\tilde{\pi}_1([\under{a}],x))\!\in\!\partial M, \,
v\!\in\!\ga_{([\under{a}],x)}, \, w\!\in\!\ga_{([\under{a}],x)}^{\perp}. \notag
\end{gather}
Since the bundle map
$$T\PP\!=\!\ga\!\oplus\!\ga^{\perp} \lra \ga^{\otimes2}\!\oplus\!\ga^{\perp}, \qquad
(v,w)\lra\big(v^{\otimes2},w),$$
is two-to-one, outside of the proper subbundle $\ga^{\perp}$,
\begin{gather}\label{cusp2_e4}
{\cal C}_{{\cal Z}_1}(\vph)=2\cdot N(\al),
\end{gather}
by \e_ref{cusp2_e2}, \e_ref{cusp2_e3}, and a rescaling and cobordism argument
as in Subsection~3.1 of~\cite{Z1}.
Suppose next that $([\under{a}],x)\!\in\!{\cal T}_1$.
Let $N_{(\under{a},x)}$ be the normal bundle of ${\cal T}_1$ in
${\cal K}_1'$ at $([\under{a}],x)$.
Then, with appropriate identifications, for some $\be_2,\be_4\!\in\!\Bbb{C}^*$
and $C\!\in\!\Bbb{R}^+$,
\begin{gather}\label{cusp2_e5}
\big|\vph([\under{a}],x;u,v,w)-\al_0(u,v,w)\big| \le C\big(|v|^5\!+\!|w|^2)
\quad\forall\, u\!\in\!N_{(\under{a},x)},\,
v\!\in\!\ga_{([\under{a}],x)},\, w\!\in\!\ga_{([\under{a}],x)}^{\perp},\\
\hbox{where}\quad
\al_0(u,v,w)=\big(\frac{1}{6}uv^3+\frac{1}{4}\be_4v^4,
\frac{1}{2}uv^2+\be_4v^3,\be_2w).
\notag
\end{gather}
Since the polynomial~$\al_0$ is four-to-one near the origin,
it follows from \e_ref{cusp2_e5} that each point of ${\cal Z}_0\!\approx\!{\cal T}_1$
contributes~$4$ to ${\cal C}_{{\cal Z}_0}(\vph)$.
From~\e_ref{cusp2_e4} and Lemmas~\ref{tacnode1_lmm} and~\ref{cusp2_lmm2}, 
we conclude that
\begin{equation}\label{cusp2_e6}
{\cal C}_{\partial M}(\vph)
={\cal C}_{{\cal Z}_1}(\vph)+{\cal C}_{{\cal Z}_0}(\vph)
=2\cdot 152+4|{\cal T}_1|=2\cdot 152+ 4\cdot 200=1104.
\end{equation}
The lemma follows from \e_ref{cusp2_e1} and~\e_ref{cusp2_e6}.

\begin{lmm}
\label{cusp2_lmm2}
If ${\cal K}_1'\!\subset\!\Bbb{P}T\PP|_{{\cal N}_2''}$ is
as in the proof of Lemma~\ref{tacnode1_lmm} and 
\begin{gather*}
\al\in\Ga\big({\cal K}_1';\Hom(\ga^{\otimes2}\!\oplus\!\ga^{\perp},
\ga_0^*\!\otimes\!\ga_1^{*\otimes4}\oplus
\ga_0^*\!\otimes\!\ga_1^{*\otimes4}\!\otimes\!T^*\PP)\big),\\
\al(v,w)=\big(0,\frac{1}{2}D_{\under{a},x}^3v,\tilde{H}_{\under{a},x}w\big)
\in \ga_0^*\!\otimes\!\ga_1^{*\otimes4}\oplus
\ga_0^*\!\otimes\!\ga_1^{*\otimes4}\!\otimes\!\ga^*
\oplus\ga_0^*\!\otimes\!\ga_1^{*\otimes4}\!\otimes\!\ga^{\perp*},
\end{gather*} 
then $N(\al)\!=\!152$.
\end{lmm}

\noindent
{\it Proof:}
Since the linear map
$$\al\!: \ga^{\perp}\lra \ga_0^*\!\otimes\!\ga_1^{*\otimes4}\!\otimes\!\ga^{\perp*}$$
is an isomorphism over ${\cal K}_1'$,
\begin{gather}\label{cusp2_e11}
N(\al)=N(\tilde{\al}), \qquad\hbox{where}\\
\tilde{\al}\in\Ga\big({\cal K}_1';\Hom(\ga^{\otimes2},
\ga_0^*\!\otimes\!\ga_1^{*\otimes4}\oplus
\ga_0^*\!\otimes\!\ga_1^{*\otimes4}\!\otimes\!\ga^*)\big), \quad
\al\big([\under{a}],x;v,w)=\big(0,\frac{1}{2}D_{\under{a},x}^3v\big).\notag
\end{gather}
Similarly to the proof of the Lemma~\ref{node2_lmm2},
\begin{equation}\label{cusp2_e12}\begin{split}
N(\tilde{\al})&=\blr{c\big(\ga_0^*\!\otimes\!\ga_1^{*\otimes4}\oplus
\ga_0^*\!\otimes\!\ga_1^{*\otimes4}\!\otimes\!\ga^*\big)
c(\ga^{\otimes2})^{-1},{\cal K}_1'}
-{\cal C}_{\tilde{\al}^{-1}(0)}(\tilde{\al}^{\perp})\\
&=\lr{3\la^3\!+\!(8y\!+\!23a)\la^2\!+\!(7y^2\!+\!41ya\!+\!61a^2)\la,
\Bbb{P}T\PP|_{{\cal N}_1''}}-{\cal C}_{{\cal T}_1}(\tilde{\al}^{\perp})\\
&=\lr{7y^2+17ya+10a^2,{\cal N}_1''}-{\cal C}_{{\cal T}_1}(\tilde{\al}^{\perp})\\
&=\big(7|{\cal N}_1|+17|{\cal N}_{1,1}|+10\lr{a,{\cal N}_{1,1}'})
-{\cal C}_{{\cal T}_1}(\tilde{\al}^{\perp})
=7\cdot 27+17\cdot 9 +10\cdot1 -{\cal C}_{{\cal T}_1}(\tilde{\al}^{\perp}).
\end{split}\end{equation}
Since the section $D^3$ is transverse to the zero set, so is
the section~$\tilde{\al}^{\perp}$ if $\nu$ is generic.
Thus,
\begin{equation}\label{cusp2_e15}
{\cal C}_{\tilde{\al}^{-1}(0)}(\tilde{\al}^{\perp})= \,
^{\pm}\big|\tilde{\al}^{-1}(0)\big|=|{\cal T}_1|.
\end{equation}
The lemma follows from \e_ref{cusp2_e11}-\e_ref{cusp2_e15}
along with Lemma~\ref{tacnode1_lmm}.

\subsection{Quartics with Three Simple Nodes}
\label{p3_subs}

\noindent
In this subsection we compute the last number of Table~\ref{quartics_table}.
We start with the following structural lemma.

\begin{lmm}
\label{node3_lmm0}
Let ${\cal N}_1''\!\subset\!{\cal D}\!\times\!\PP_1$ be 
as in the proof of Lemma~\ref{tacnode1_lmm} and let
\begin{gather*}
\tilde{\cal N}_{2;0}'=\big\{([\under{a}],x_1,x_2)\!\in\!{\cal N}_1''\!\times\!\PP_2\!:
x_1\!\neq\!x_2,\, \vph_2([\under{a}],x_1,x_2)\!=\!0\},
\qquad\hbox{where}\\
\vph_2\!\in\!\Ga\big({\cal N}_1''\!\times\!\PP_2;\ga_0^*\!\otimes\!\ga_2^{*\otimes4}\oplus
\ga_0^*\!\otimes\!\ga_2^{*\otimes4}\!\otimes\!T^*\PP_2\big), \quad
\vph_2([\under{a}],x_1,x_2)=\big(s_{\under{a}}(x_2),ds_{\under{a}}|_{x_2}\big).
\end{gather*}
If $\tilde{\cal N}_2'$ is the closure of $\tilde{\cal N}_{2;0}'$ 
in ${\cal N}_1''\!\times\!\PP_2$, then
$$\partial\tilde{\cal N}_2'\equiv\tilde{\cal N}_2'-\tilde{\cal N}_{2;0}'
=\big\{([\under{a}],x,x)\!\in\!{\cal N}_1''\!\times\!\PP_2\!:
([\under{a}],x)\!\in\!{\cal T}_1\big\}.$$
\end{lmm}

\noindent
{\it Proof:}
We will only show that if $([\under{a}],x_1,x_2)\!\in\!\partial\tilde{\cal N}_2'$,
then $x_1\!=\!x_2$ and $([\under{a}],x_1)\!\in\!{\cal T}_1$.
The converse follows from the proofs of Lemmas~\ref{node3_lmm1} and~\ref{node3_lmm2}.
Suppose $([\under{a}],x_1,x_2)\!\in\!\partial\tilde{\cal N}_2'$.
Since the section $\vph_2$ is continuous, $x_2\!=\!x_1$ by definition of 
$\partial\tilde{\cal N}_2'$.
If $([\under{a}],x_1)\!\in\!{\cal N}_1''\!-\!{\cal K}_1'$, by 
\e_ref{node2_e2} and~\e_ref{node2_e3} and with appropriate identifications,
$$|\vph_2([\under{a}],x_1,v)|\ge C([\under{a}],x_1)^{-1}|v|$$
for all $v\!\in\!T_{x_1}\PP_1$ sufficiently small.
Thus, $([\under{a}],x_1,x_1)$ is not in the closure of~$\tilde{\cal N}_{2;0}$.
Suppose next that $(\under{a},x_1)\!\in\!{\cal K}_1'\!-\!{\cal T}_1$.
Then, by \e_ref{node2_e5}, 
\begin{equation}\label{node3_e0}
|\vph_2([\under{a}],x_1,u,v,w)|\ge C([\under{a}],x_1)^{-1}(|v|^3\!+\!|w|)
\end{equation}
for all $u\!\in\!N_{(\under{a},x_1)}$, $v\!\in\!{\cal L}_{(\under{a},x_1)}$,
and $\!\in\!{\cal L}_{(\under{a},x_1)}^{\perp}$ sufficiently small.
In this case, $N$ is the normal bundle of ${\cal K}_1'$, viewed as a submanifold
of~${\cal N}_1''$, in~${\cal N}_1''$, while the line bundles ${\cal L}$
and ${\cal L}^{\perp}$ over ${\cal K}_1'$ are defined as in~(1)
of the proof of Lemma~\ref{node2_lmm1}.
From~\e_ref{node3_e0},
we conclude that $([\under{a}],x_1,x_1)$ is not in the closure of~$\tilde{\cal N}_{2;0}$.

\begin{lmm}
\label{node3_lmm1}
The number $|{\cal N}_3|$ of plane quartics that have three nodes and pass
through $11$ points in general position is~$675$.
\end{lmm}

\noindent
{\it Proof:}
With notation as in the statement of Lemma~\ref{node3_lmm0}, let
\begin{gather*}
M\!=\!\tilde{\cal N}_2'\!\times\!\PP_3, ~~
M^0\!=\!\big\{([\under{a}],x_1,x_2,x_3)\!\in\!M\!:x_3\!\neq\!x_1,x_2\},
~~ \partial M\!=\!M\!-\!M^0, ~~
\tilde{\cal N}_3\!=\!\vph_3^{-1}(0)\cap M^0,\\
\hbox{where}~~~
\vph_3\!\in\!\Ga\big(M;\ga_0^*\!\otimes\!\ga_3^{*\otimes4}\oplus
\ga_0^*\!\otimes\!\ga_3^{*\otimes4}\!\otimes\!T^*\PP_3\big),  ~~
\vph_3([\under{a}],x_1,x_2,x_3)=\big(s_{\under{a}}(x_3),ds_{\under{a}}|_{x_3}\big),
~~ \ga_3\!=\!\pi_3^*\ga_{\PP_3},
\end{gather*}
and $\pi_3\!:M\!\lra\!\PP_3$ is the projection onto the last component.
Since $\vph_3|_{M^0}$ is transverse to the zero set,
\begin{equation}\label{node3_e1}\begin{split}
|\tilde{\cal N}_3| &=\, ^{\pm}\big|\vph_3^{-1}(0)\cap M^0\big|
= \blr{e(\ga_0^*\!\otimes\!\ga_3^{*\otimes4}\oplus
\ga_0^*\!\otimes\!\ga_3^{*\otimes4}\!\otimes\!T^*\PP_3),M}
-{\cal C}_{\partial M}(\vph_3)\\
&= \blr{(y\!+\!4a_3)(y^2\!+\!5ya_3\!+\!7a_3^2),\tilde{\cal N}_2'\!\times\!\PP_3}
-{\cal C}_{\partial M}(\vph_3)\\
&= 27\lr{y,\tilde{\cal N}_2'}-{\cal C}_{\partial M}(\vph_3)
=27\cdot2|{\cal N}_2|-{\cal C}_{\partial M}(\vph_3)
=27\cdot 450-{\cal C}_{\partial M}(\vph_3),
\end{split}\end{equation}
where $a_3\!=\!\pi_2^*c_1(\ga_{\PP_3}^*)$.
In order to determine ${\cal C}_{\partial M}(\vph_3)$,
we split $\partial M$ into five strata:
\begin{alignat*}{2}
&{\cal Z}_{1,i}=\big\{([\under{a}],x_1,x_2,x_3)\!: x_3\!=\!x_i, \, x_3\!\neq\!x_j, \,
([\under{a}],x_i)\!\in\!{\cal N}_1''\!-\!{\cal K}_1'\big\},
&\quad& \{i,j\}=\{1,2\};\\
&{\cal Z}_{0,i}=\big\{([\under{a}],x_1,x_2,x_3)\!: x_3\!=\!x_i, \, x_3\!\neq\!x_j, \,
([\under{a}],x_i)\!\in\!{\cal K}_1'\!-\!{\cal T}_1\big\}\approx{\cal K}_2,
&\quad& \{i,j\}=\{1,2\};\\
&{\cal Z}_{0,12}=\big\{([\under{a}],x,x,x)\!:([\under{a}],x)\!\in\!{\cal T}_1\big\}. &&
\end{alignat*}
Note that Lemma~\ref{node3_lmm0} implies that the union of these five
spaces is indeed~$\partial M$.
Similarly to the proof of Lemma~\ref{node2_lmm1}, we have 
\begin{gather}\label{node3_e2}
{\cal C}_{{\cal Z}_{1,1}}(\vph_3)=
{\cal C}_{{\cal Z}_{1,2}}(\vph_3)=N(\al),   \qquad\hbox{where}\\
\al\in\Ga\big(\tilde{\cal N}_2';\Hom(T\PP_1,\ga_0^*\!\otimes\!\ga_1^{*\otimes4}\oplus
\ga_0^*\!\otimes\!\ga_1^{*\otimes4}\!\otimes\!T^*\PP_1)\big),
\quad \al([\under{a}],x_1,x_2;v)=(0,H_{\under{a},x_1}v),  \notag
\end{gather}
while
\begin{equation}\label{node3_e3}
{\cal C}_{{\cal Z}_{0,1}}(\vph_3)=
{\cal C}_{{\cal Z}_{0,2}}(\vph_3)= 3|{\cal K}_2|.
\end{equation}
Finally, suppose that $([\under{a}],x)\!\in\!{\cal T}_1$.
Let $N_{(\under{a},x)}^1$ and $N_{(\under{a},x)}^2$ be the normal bundles of 
${\cal T}_1$ in ${\cal K}_1'$ and of ${\cal K}_1'$ in ${\cal N}_1''$,
respectively, at $([\under{a}],x)$.
Let ${\cal L}_{(\under{a},x)}$ and ${\cal L}_{(\under{a},x)}^{\perp}$
be as in the proof of Lemma~\ref{node2_lmm1}.
Then, with appropriate identifications, for some $\be_2,\be_4\!\in\!\Bbb{C}^*$,
$C\!\in\!\Bbb{R}^+$, and $i\!=\!2,3$,
\begin{gather}\label{node3_e4}
\big|\vph_i([\under{a}],x;u_1,u_2,v_i,w_i)-\al_0(u_1,u_2,v_i,w_i)\big| 
\le C\big(|v_i|^5\!+\!|w_i|^2)\\
\forall~ u_1\!\in\!N_{(\under{a},x)}^1,\, u_2\!\in\!N_{(\under{a},x)}^2, \,
v_i\!\in\!{\cal L}_{(\under{a},x)},\, w_i\!\in\!{\cal L}_{(\under{a},x)}^{\perp}, \notag\\
\hbox{where}\quad
\al_0(u_1,u_2,v,w)=\big(\frac{1}{6}u_1v^3+\frac{1}{2}u_2v^2+\frac{1}{12}\be_4v^4,
\frac{1}{2}u_1v^2+u_2v+\frac{1}{3}\be_4v^3,\be_2w).
\notag
\end{gather}
Since $\tilde{\cal N}_{2;0}'\!=\!\vph_2^{-1}(0)$, the $\vph_3$-contribution of 
$([\under{a}],x,x,x)$ is the number of small solutions of the system
\begin{equation}\label{node3_e5}
\begin{cases}
\vph_2(u_1,u_2,v_2,w_2)=0\\
\vph_3(u_1,u_2,v_3,w_3)=t\, \nu(u_1,u_2,v_2,w_2,v_3,w_3)
\end{cases} \qquad
\begin{aligned}
&(u_1,u_2,v_2,w_2,v_3,w_3)\in\Bbb{C}^6\\
& (v_2,w_2)\neq(0,0),
\end{aligned}
\end{equation}
for a generic $\nu\!\in\!\Bbb{C}^3$ and $t\!\in\!\Bbb{R}^+$ sufficiently small.
By \e_ref{node3_e4} and
a rescaling and cobordism argument as in Subsection~3.1 of~\cite{Z1},
the number of small solutions of~\e_ref{node3_e5}
is the same as the number of solutions of the system
\begin{equation}\label{node3_e6}
\begin{cases}
\frac{1}{6}u_1v_2^3+\frac{1}{2}u_2v_2^2+\frac{1}{12}\be_4v_2^4=0\\
\frac{1}{2}u_1v_2^3+u_2v_2^2+\frac{1}{3}\be_4v_2^4=0\\
\frac{1}{6}u_1v_3^3+\frac{1}{2}u_2v_3^2+\frac{1}{12}\be_4v_3^4=\nu\\
\frac{1}{2}u_1v_3^3+u_2v_3^2+\frac{1}{3}\be_4v_3^4=0\\
\end{cases} \qquad
(u_1,u_2,v_2,v_3)\in\Bbb{C}^{*4},
\end{equation}
for a generic $\nu\!\in\!\Bbb{C}$.
Dividing the first two equations by $v_2^2$ and the last equation by $v_3^2$
and then solving for $u_2$ and $u_1$ in terms of $v_1$ and $v_2$,
we find that the system~\e_ref{node3_e6} is equivalent~to
\begin{equation}\label{node3_e6b}
\begin{cases}
u_1=-\be_4v_2\\
u_2=\frac{1}{6}\be_4v_2^2\\
v_2=v_3~\hbox{or}~v_2=2v_3\\
-\frac{1}{6}v_2v_3^3+\frac{1}{12}v_2^2v_3^2+\frac{1}{12}v_4^4=\nu
\end{cases} \qquad
(u_1,u_2,v_2,v_3)\in\Bbb{C}^4.
\end{equation}
If $v_2\!=\!v_3$, the last equation has no solutions for $\nu\!\neq\!0$.
On the other hand, if $v_2\!=\!2v_3$, the last equation in~\e_ref{node3_e6b}
has four solutions.
We conclude that
\begin{equation}\label{node3_e7}
{\cal C}_{{\cal Z}_{0,12}}(\vph_3)=4|{\cal T}_1|.
\end{equation}
From \e_ref{node3_e2}, \e_ref{node3_e3}, and \e_ref{node3_e7},
along with Lemmas~\ref{tacnode1_lmm} and \ref{cusp2_lmm1}, we conclude~that
\begin{equation}\label{node3_e8}
{\cal C}_{\partial M}(\vph_3)
=2\, {\cal C}_{{\cal Z}_{1,1}}(\vph_3)+2\, {\cal C}_{{\cal Z}_{0,1}}(\vph_3)
+{\cal C}_{{\cal Z}_{0,12}}(\vph_3)
=2\cdot 1130 + 6\cdot 840 + 4\cdot 200=  8100.
\end{equation}
The lemma follows from \e_ref{node3_e1} and \e_ref{node3_e8},
since ${\cal N}_3\!=\!\tilde{\cal N}_3/S_3$.

\begin{lmm}
\label{node3_lmm2}
If $\tilde{\cal N}_2'\!\subset\!{\cal D}\!\times\!\PP_1\!\times\!\PP_2$ is
as in Lemma~\ref{node3_lmm0} and 
$$\al\in\Ga\big(\tilde{\cal N}_2';\Hom(T\PP_1,\ga_0^*\!\otimes\!\ga_1^{*\otimes4}\oplus
\ga_0^*\!\otimes\!\ga_1^{*\otimes4}\!\otimes\!T^*\PP_1)\big),
\quad \al([\under{a}],x_1,x_2;v)=(0,H_{\under{a},x_1}v),$$ 
then $N(\al)\!=\!1130$.
\end{lmm}

\noindent
{\it Proof:} We put 
$$M=\Bbb{P}T\PP_1|_{\tilde{\cal N}_2'}, \quad
\partial M=\big\{([\under{a}],x_1,x_2)\!\in\!M\!: \tilde{H}_{\under{a},x}\!=\!0\big\},$$
where $\tilde{H}_{\cdot,\cdot}$ is as in the proof of Lemma~\ref{tacnode1_lmm}.
Using Lemma~\ref{node3_lmm0}, we split $\partial M$ into two subsets:
\begin{alignat*}{1}
&{\cal Z}_{0,1}=\big\{([\under{a}],x_1,x_2)\!\in\!M\!: \tilde{\pi}_1([\under{a}],x_1)\!\neq\!x_2,\,
([\under{a}],\under{x}_1)\!\in\!{\cal K}_1'\!-\!{\cal T}_1\big\},\\
&{\cal Z}_{0,2}=\big\{([\under{a}],x_1,x_2)\!\in\!M\!: \tilde{\pi}_1([\under{a}],x_1)\!=\!x_2,\,
([\under{a}],\under{x}_1)\!\in\!{\cal T}_1\big\},
\end{alignat*}
where $\tilde{\pi}_1$ is as in the proof of Lemma~\ref{cusp2_lmm1}.
Here ${\cal K}_1'$ and ${\cal T}_1$ are viewed as subspaces of $\Bbb{P}T\PP|_{{\cal N}_1''}$,
as defined in the proof of Lemma~\ref{tacnode1_lmm}.
Let
$$\tilde{\al}=(0,\tilde{H})
\in \Ga\big(M;\Hom(\ga,\ga_0^*\!\otimes\!\ga_1^{*\otimes4}\oplus
\ga_0^*\!\otimes\!\ga_1^{*\otimes4}\!\otimes\!T^*\PP_1)\big)$$
be the section induced by $\al$.
Similarly to the proof of Lemma~\ref{node2_lmm2},
\begin{equation}\label{node3_e21}\begin{split}
N(\al)&=\blr{c\big(\ga_0^*\!\otimes\!\ga_1^{*\otimes4}\oplus
\ga_0^*\!\otimes\!\ga_1^{*\otimes4}\!\otimes\!T^*\PP_1\big)c(T\PP_1)^{-1},\tilde{\cal N}_2'}
-{\cal C}_{\tilde{\al}^{-1}(0)}(\tilde{\al}^{\perp})\\
&=\lr{3y+6a_1,\tilde{\cal N}_2'}-{\cal C}_{\partial M}(\tilde{\al}^{\perp})
=\big(3\cdot2|{\cal N}_2|+6|{\cal N}_{2,1}|)-{\cal C}_{\partial M}(\tilde{\al}^{\perp}).
\end{split}\end{equation}
As in the proof of Lemma~\ref{node2_lmm2}, we have
\begin{equation}\label{node3_e22}
{\cal C}_{{\cal Z}_{0,1}}(\tilde{\al}^{\perp})= \,
^{\pm}\big|{\cal Z}_{0,1}\big|=|{\cal K}_2|.
\end{equation}
On the other hand, suppose $([\under{a}],x_1,x_2)\!\in\!{\cal Z}_{0,2}$
and thus $([\under{a}],x_1)\!\in\!{\cal T}_1$, while
$x_2\!=\!\tilde{\pi}_1([\under{a}],x_1)$.
Then, with identifications similar to the ones used at the end of the proof of
Lemma~\ref{node3_lmm1}, the $\tilde{\al}^{\perp}$-contribution of $([\under{a}],x_1,x_2)$
is the number of small solutions of the system
\begin{equation}\label{node3_e23}
\begin{cases}
\vph_2(u_1,u_2,v_2,w_2)=0\\
\tilde{\al}^{\perp}(u_1,u_2,w_3)=t\, \nu(u_1,u_2,v_2,w_2,w_3)
\end{cases} \qquad
\begin{aligned}
&(u_1,u_2,v_2,w_2,w_3)\in\Bbb{C}^5\\
& (v_2,w_2)\neq(0,0),
\end{aligned}
\end{equation}
for a generic $\nu\!\in\!\Bbb{C}^2$ and $t\!\in\!\Bbb{R}^+$ sufficiently small.
In this case, $w_3\!\in\!\ga_{(\under{a},x_1)}^*\!\otimes\!\ga_{(\under{a},x_1)}^{\perp}$.
For a good choice of identifications
\begin{equation}\label{node3_e24}
\tilde{\al}^{\perp}(u_1,u_2,v_2,w_2,w_3)=(u_2,w_3).
\end{equation}
By the $i\!=\!2$ case of \e_ref{node3_e4} and \e_ref{node3_e24},
the number of small solutions of the system~\e_ref{node3_e23} is 
the same as the number of solutions of the system
$$\begin{cases}
\frac{1}{6}u_1v_2^3+\frac{1}{2}u_2v_2^2+\frac{1}{12}\be_4v_2^4=0\\
\frac{1}{2}u_1v_2^3+u_2v_2^2+\frac{1}{3}\be_4v_2^4=0\\
u_2=\nu
\end{cases} \qquad
(u_1,u_2,v_2)\in\Bbb{C}^{*3}, $$
for a generic $\nu\!\in\!\Bbb{C}$.
Thus, each point of ${\cal Z}_{0,2}$ contributes two, and
\begin{equation}\label{node3_e26}
{\cal C}_{{\cal Z}_{0,2}}(\tilde{\al}^{\perp})=2|{\cal T}_1|.
\end{equation}
The lemma follows from \e_ref{node3_e21}, \e_ref{node3_e22}, and \e_ref{node3_e26},
along with Lemmas~\ref{tacnode1_lmm}, \ref{node2_lmm1}, and~\ref{cusp2_lmm1}.

\subsection{Generalization to Arbitrary-Degree Curves}
\label{curvgen_subs}

\noindent
The computations in the previous subsections generalize to higher-degree
curves, as well as to other types of singularities.
We list the results of the generalization to arbitrary-degree curves 
in Table~\ref{curvegen_table}.
The number in the third column is the lowest value of the degree~$d$ for which
the formula given in the last column is applicable.
Note that in the cases when this number is higher than one,
the constraints are $-1$ points for $d\!=\!1$ and two points for $d\!=\!2$.
So, the corresponding count of curves makes no sense for $d\!=\!1$,
while for $d\!=\!2$ this is a count of structures on the double line through two
distinct points in~$\PP$.
The number in the fourth column is the difference between
$$\dim(d) \equiv \dim\big\{\hbox{deg.}-\!d~\hbox{curves}\big\} = \frac{d(d\!+\!3)}{2}$$
and the number of points in general position.
Below we state the changes that are needed to be made in the above lemmas to obtain
these results.

\begin{table}
\begin{center}
\begin{tabular}{||c|c|c|c|c||}
\hline\hline
set& singularities& $d\!\ge\!$& co-\# pts& cardinality\\
\hline\hline
${\cal N}_1$& 1 node& 1& 1& $3(d\!-\!1)^2$\\
\hline
${\cal N}_{1,1}$& 1 node on a fixed line& 1& 2& $3(d\!-\!1)$\\
\hline
${\cal K}_1$& 1 cusp& 1& 2& $12(d\!-\!1)(d\!-\!2)$\\
\hline
${\cal K}_{1,1}$& 1 cusp on a fixed line& 3& 3& $4(2d\!-\!3)$\\
\hline
${\cal T}_1$& 1 tacnode& 3& 3& $2(25d^2\!-\!96d\!+\!84)$\\
\hline\hline
${\cal N}_2$& 2 nodes& 1& 2& $3(d\!-\!1)(d\!-\!2)(3d^2\!-\!3d\!-\!11)/2$\\
\hline
${\cal N}_{2,1}$& 2 nodes, one on a fixed line& 3& 3& $9d^3\!-\!27d^2\!-\!d\!+\!30$\\
\hline
${\cal K}_2$& 1 node and 1 cusp& 3& 3& $12(d\!-\!3)(3d^3\!-\!6d^2\!-\!11d\!+\!18)$\\
\hline\hline
${\cal N}_3$& 3 nodes& 3& 3& 
\begin{small}$\!\!\!(9d^6\!-\!54d^5\!+\!9d^4\!+\!423d^3\!-\!458d^2
\!-\!829d\!+\!1050)/2\!\!\!$\end{small} \\
\hline\hline
\end{tabular}
\caption{Some Characteristic Numbers of Degree-$d$ Plane Curves}
\label{curvegen_table}
\end{center}
\end{table}

\subsubsection{The Numbers ${\cal N}_1$ and ${\cal N}_{1,1}$}
\label{n1gen_subsub}

\noindent
In order to compute the number ${\cal N}_1$, we take ${\cal D}\!\approx\!\Bbb{P}^1$ 
to be the subspace of degree-$d$ plane curves that pass through a set of $\dim(d)\!-\!1$ 
points in general position.
We define ${\cal N}_1$ as in~(1) of the proof of Lemma~\ref{node1_lmm}, except now
$$\vph\in\Ga\big({\cal D}\!\times\!\PP;
\ga_0^*\!\otimes\!\ga_1^{*\otimes d}\oplus
\ga_0^*\!\otimes\!\ga_1^{*\otimes d}\!\otimes\!T^*\PP\big).$$
Since $\vph$ is transverse to the zero set, we obtain
\begin{equation*}\begin{split}
|{\cal N}_1|=|\vph^{-1}(0)|&
=\blr{e\big(\ga_0^*\!\otimes\!\ga_1^{*\otimes d}\oplus
\ga_0^*\!\otimes\!\ga_1^{*\otimes d}\!\otimes\!T^*\PP\big),
{\cal D}\!\times\!\PP}\\
&=\blr{(y\!+\!da)(y^2\!+\!(2d\!-\!3)ya\!+\!(d^2\!-\!3d\!+\!3)a^2),{\cal D}\!\times\!\PP}
=3(d\!-\!1)^2.
\end{split}\end{equation*}
With the analogous changes in (2) of the proof of Lemma~\ref{node1_lmm},
we find that
\begin{equation*}\begin{split}
|{\cal N}_{1,1}|=|\vph^{-1}(0)|&
=\blr{e\big(\ga_0^*\!\otimes\!\ga_1^{*\otimes d}\oplus
\ga_0^*\!\otimes\!\ga_1^{*\otimes d}\!\otimes\!T^*\PP\big),
{\cal D}\!\times\!\Bbb{P}^1}\\
&=\blr{(y\!+\!da)(y^2\!+\!(2d\!-\!3)ya\!+\!(d^2\!-\!3d\!+\!3)a^2),
{\cal D}\!\times\!\Bbb{P}^1}=3(d\!-\!1).
\end{split}\end{equation*}

\subsubsection{The Numbers ${\cal K}_1$ and ${\cal K}_{1,1}$}
\label{k1gen_subsub}

\noindent
We take ${\cal D}\!\approx\!\Bbb{P}^2$ to be the subspace of degree-$d$ plane curves 
that pass through a set of $\dim(d)\!-\!2$ points in general position.
We define ${\cal N}_1'$ and ${\cal K}_1$ as in~(1) of the proof of Lemma~\ref{cusp1_lmm}, 
except now
$$H_{\under{a},x}\in \Ga\big({\cal N}_1';
\Hom(T\PP,\ga_0^*\!\otimes\!\ga_1^{*\otimes d}\!\otimes\!T^*\PP)\big)
\quad\hbox{and}\quad
\vph\in\Ga\big({\cal N}_1';
(\ga_0^*\!\otimes\!\ga_1^{*\otimes d}\!\otimes\!\La^2T^*\PP)^{\otimes2}\big).$$
Since $\vph$ is transverse to the zero set,
\begin{equation*}\begin{split}
|{\cal K}_1|=|\vph^{-1}(0)|&
=\blr{e\big((\ga_0^*\!\otimes\!\ga_1^{*\otimes d}\!\otimes\!\La^2T^*\PP)^{\otimes2}\big),
{\cal N}_1'}\\
&=2\blr{y\!+\!(d\!-\!3)a,{\cal N}_1'}
=2\, \big(|{\cal N}_1|\!+\!(d\!-\!3)|{\cal N}_{1,1}|\big)\\
&=2\big(3(d\!-\!1)^2+(d\!-\!3)\cdot3(d\!-\!1)\big) = 12(d\!-\!1)(d\!-\!2).
\end{split}\end{equation*}
With the analogous changes in (2) of the proof of Lemma~\ref{cusp1_lmm},
we find that
\begin{equation*}\begin{split}
|{\cal K}_{1,1}|&
=\blr{e\big((\ga_0^*\!\otimes\!\ga_1^{*\otimes d}\!\otimes\!\La^2T^*\PP)^{\otimes2}\big),
{\cal N}_{1,1}'}\\
&=2\blr{y\!+\!(d\!-\!3)a,{\cal N}_{1,1}'}
=2\, \big(|{\cal N}_{1,1}|\!+\!(d\!-\!3)\lr{a,{\cal N}_{1,1}'}\big)\\
&=2\big(3(d\!-\!1)+(d\!-\!3)\big)=4(2d\!-\!3).
\end{split}\end{equation*}

\subsubsection{The Number ${\cal T}_1$}
\label{t1gen_subsub}

\noindent
In this case, we take ${\cal D}\!\approx\!\Bbb{P}^3$ to be the subspace of degree-$d$ 
plane curves  that pass through a set of $\dim(d)\!-\!3$ points in general position.
We define ${\cal N}_1''$, $M$, ${\cal K}_1'$, and ${\cal T}_1$ 
as in the proof of Lemma~\ref{tacnode1_lmm}, except now
$$\tilde{H}_{\cdot,\cdot}\in 
\Ga\big(M;\Hom(\ga,\ga_0^*\!\otimes\!\ga_1^{*\otimes d}\!\otimes\!T^*\PP)\big)
\quad\hbox{and}\quad
\vph\in\Ga\big(M;\Hom(\ga^{\otimes3},\ga_0^*\!\otimes\!\ga_1^{*\otimes d})\big).$$
Since the sections $\vph$ and $\tilde{H}_{\cdot,\cdot}$ are transverse 
to the zero set, we obtain
\begin{equation*}\begin{split}
|{\cal T}_1|=|\vph^{-1}(0)|
&=\blr{e\big(\ga^{*\otimes3}\!\otimes\!\ga_0^*\!\otimes\!\ga_1^{*\otimes d}\big)
e\big(\ga^*\!\otimes\!\ga_0^*\!\otimes\!\ga_1^{*\otimes d}\!\otimes\!T^*\PP\big),M}\\
&=\blr{3\la^3+\big(7y\!+\!(7d\!-\!9)a\big)\la^2+
\big(5y^2\!+\!(10d\!-\!12)ya\!+\!(5d^2\!-\!12d\!+\!9)a^2\big)\la,M}\\
&=\blr{5y^2+(10d\!-\!33)ya+(5d^2\!-\!33d\!+\!54)a^2,{\cal N}_1''}\\
&=5|{\cal N}_1|+(10d\!-\!33)|{\cal N}_{1,1}|+(5d^2\!-\!33d\!+\!54)\lr{a,{\cal N}_{1,1}'}\\
&=5\cdot3(d\!-\!1)^2+(10d\!-\!33)\cdot3(d\!-\!1)+(5d^2\!-\!33d\!+\!54)
=2\big(25d^2\!-\!96d\!+\!84).
\end{split}\end{equation*}

\subsubsection{The Numbers ${\cal N}_2$ and ${\cal N}_{2,1}$}
\label{n2gen_subsub}

\noindent
In order to compute the number ${\cal N}_2$, we take ${\cal D}\!\approx\!\Bbb{P}^2$ 
to be the subspace of degree-$d$ plane curves that pass through a set of $\dim(d)\!-\!2$ 
points in general position.
We define ${\cal N}_1'$,  $M$, $\partial M$, $\tilde{\cal N}_2$, ${\cal Z}_1$,
${\cal Z}_0$, and $\al$ as in (1) of  the proof of Lemma~\ref{node2_lmm1}, except now
$$\vph\!\in\!\Ga(M;\ga_0^*\!\otimes\!\ga_2^{*\otimes d}\oplus
\ga_0^*\!\otimes\!\ga_2^{*\otimes d}\!\otimes\!T^*\PP_2), \qquad
\al\in\Ga\big({\cal N}_1';\Hom(T\PP,\ga_0^*\!\otimes\!\ga_1^{*\otimes d}\oplus
\ga_0^*\!\otimes\!\ga_1^{*\otimes d}\!\otimes\!T^*\PP)\big).$$
Since $\vph|_{M^0}$ is transverse to the zero set,
\begin{equation}\label{node2gen_e1}\begin{split}
|\tilde{\cal N}_2| &=\, ^{\pm}\big|\vph^{-1}(0)\!\cap\!M^0\big|
= \blr{e(\ga_0^*\!\otimes\!\ga_2^{*\otimes d}\oplus
\ga_0^*\!\otimes\!\ga_2^{*\otimes d}\!\otimes\!T^*\PP_2),M}
-{\cal C}_{\partial M}(\vph)\\
&= \blr{(y\!+\!da_2)(y^2\!+\!(2d\!-\!3)ya_2\!+\!(d^2\!-\!3d\!+\!3)a_2^2),
{\cal N}_1'\!\times\!\PP_2}
-{\cal C}_{\partial M}(\vph)\\
&= 3(d\!-\!1)^2\lr{y,{\cal N}_1'}-{\cal C}_{\partial M}(\vph)
=3(d\!-\!1)^2|{\cal N}_1|-
\big({\cal C}_{{\cal Z}_0}(\vph)\!+\!{\cal C}_{{\cal Z}_1}(\vph)\big).
\end{split}\end{equation}
As in (1) of the proof of Lemma~\ref{node2_lmm1}, we have
\begin{equation}\label{node2gen_e2}
{\cal C}_{{\cal Z}_1}(\vph)=N(\al) \qquad\hbox{and}\qquad
{\cal C}_{{\cal Z}_0}(\vph)=3\big|{\cal K}_1\big|.
\end{equation}
Similarly to (1) of the proof of Lemma~\ref{node2_lmm2},
\begin{gather*}\begin{split}
N(\al)&=\blr{c\big(\ga_0^*\!\otimes\!\ga_1^{*\otimes d}\oplus
\ga_0^*\!\otimes\!\ga_1^{*\otimes d}\!\otimes\!T^*\PP\big)c(T\PP)^{-1},{\cal N}_1'}
-{\cal C}_{\tilde{\al}^{-1}(0)}(\tilde{\al}^{\perp})\\
&=\lr{3y+3(d\!-\!2)a,{\cal N}_1'}-{\cal C}_{\partial M}(\tilde{\al}^{\perp})
=\big(3|{\cal N}_1|+3(d\!-\!2)|{\cal N}_{1,1}|)-{\cal C}_{\partial M}(\tilde{\al}^{\perp}),
\end{split}\\
\hbox{where}\qquad {\cal C}_{\partial M}(\tilde{\al}^{\perp})=\big|{\cal K}_1\big|.
\end{gather*}
Combining these observations with \e_ref{node2gen_e1} and \e_ref{node2gen_e2},
we obtain
\begin{equation*}
|\tilde{\cal N}_2| 
= 3d(d\!-\!2)|{\cal N}_1|-3(d\!-\!2)|{\cal N}_{1,1}| -2|{\cal K}_1|
= 3(d\!-\!1)(d\!-\!2)\big(3d^2\!-\!3d\!-\!11\big).
\end{equation*}\\

\noindent
With the analogous modifications in (2) of the proof of Lemma~\ref{node2_lmm1},
we obtain
\begin{gather}
\label{node2gen_e11}\begin{split}
|{\cal N}_{2,1}| &=\, ^{\pm}\big|\vph^{-1}(0)\!\cap\! M^0\big|
= \blr{e(\ga_0^*\!\otimes\!\ga_2^{*\otimes d}\oplus
\ga_0^*\!\otimes\!\ga_2^{*\otimes d}\!\otimes\!T^*\PP_2),{\cal N}_{1,1}'\!\times\!\PP_2}
-{\cal C}_{\partial M}(\vph)\\
&= 3(d\!-\!1)^2\lr{y,{\cal N}_{1,1}'}-{\cal C}_{\partial M}(\vph)
=3(d\!-\!1)^2
|{\cal N}_{1,1}|-\big({\cal C}_{{\cal Z}_0}(\vph)\!+\!{\cal C}_{{\cal Z}_1}(\vph)\big),
\end{split}\\
\label{node2gen_e12}
\hbox{where}\qquad
{\cal C}_{{\cal Z}_1}(\vph)=N(\al) \qquad\hbox{and}\qquad
{\cal C}_{{\cal Z}_0}(\vph)=3\big|{\cal K}_{1,1}\big|.
\end{gather}
By the argument in (2) of the proof of Lemma~\ref{node2_lmm2},
\begin{gather*}\begin{split}
N(\al)&=\blr{c\big(\ga_0^*\!\otimes\!\ga_1^{*\otimes d}\oplus
\ga_0^*\!\otimes\!\ga_1^{*\otimes d}\!\otimes\!T^*\PP\big)c(T\PP)^{-1},{\cal N}_{1,1}'}
-{\cal C}_{\tilde{\al}^{-1}(0)}(\tilde{\al}^{\perp})\\
&=\lr{3y+3(d\!-\!2)a,{\cal N}_{1,1}'}-{\cal C}_{\partial M}(\tilde{\al}^{\perp})
=\big(3|{\cal N}_{1,1}|+3(d\!-\!2)\big)-{\cal C}_{\partial M}(\tilde{\al}^{\perp}),
\end{split}\\
\hbox{where}\qquad {\cal C}_{\partial M}(\tilde{\al}^{\perp})=\big|{\cal K}_{1,1}\big|.
\end{gather*}
Combining these identities with \e_ref{node2gen_e11} and \e_ref{node2gen_e12},
we obtain
\begin{equation*}
|{\cal N}_{2,1}| 
= 3d(d\!-\!2)|{\cal N}_{1,1}|-3(d\!-\!2) -2|{\cal K}_{1,1}|
= 9d^3\!-\!27d^2\!-\!d\!+\!30.
\end{equation*}

\subsubsection{The Number ${\cal K}_2$}
\label{k2gen_subsub}

\noindent
We take ${\cal D}\!\approx\!\Bbb{P}^3$ to be the subspace of degree-$d$ 
plane curves  that pass through a set of $\dim(d)\!-\!3$ points in general position.
We define ${\cal N}_1''$, ${\cal K}_1'$, $M$, $\partial M$, ${\cal K}_2$, ${\cal Z}_1$,
${\cal Z}_0$, and $\al$ as in the proof of Lemma~\ref{cusp2_lmm1}, except now
$$\vph\!\in\!\Ga(M;\ga_0^*\!\otimes\!\ga_2^{*\otimes d}\oplus
\ga_0^*\!\otimes\!\ga_2^{*\otimes d}\!\otimes\!T^*\PP_2), \quad
\al\in\Ga\big({\cal K}_1';\Hom(\ga^{\otimes2}\!\oplus\!\ga^{\perp},
       \ga_0^*\!\otimes\!\ga_1^{*\otimes d}\oplus
\ga_0^*\!\otimes\!\ga_1^{*\otimes d}\!\otimes\!T^*\PP)\big).$$
Since $\vph|_{M^0}$ is transverse to the zero set,
\begin{equation}\label{cusp2gen_e1}\begin{split}
|{\cal K}_2| &=\, ^{\pm}\big|\vph^{-1}(0)\!\cap\!M^0\big|
= \blr{e(\ga_0^*\!\otimes\!\ga_2^{*\otimes d}\oplus
\ga_0^*\!\otimes\!\ga_2^{*\otimes d}\!\otimes\!T^*\PP_2),M}
-{\cal C}_{\partial M}(\vph)\\
&= \blr{(y\!+\!da_2)(y^2\!+\!(2d\!-\!3)ya_2\!+\!(d^2\!-\!3d\!+\!3)a_2^2),
{\cal K}_1'\!\times\!\PP_2}
-{\cal C}_{\partial M}(\vph)\\
&= 3(d\!-\!1)^2\lr{y,{\cal K}_1'}-{\cal C}_{\partial M}(\vph)
=3(d\!-\!1)^2|{\cal K}_1|-
\big({\cal C}_{{\cal Z}_0}(\vph)\!+\!{\cal C}_{{\cal Z}_1}(\vph)\big).
\end{split}\end{equation}
As in (1) of the proof of Lemma~\ref{cusp2_lmm1}, we have
\begin{equation}\label{cusp2gen_e2}
{\cal C}_{{\cal Z}_1}(\vph)=2N(\al) \qquad\hbox{and}\qquad
{\cal C}_{{\cal Z}_0}(\vph)=4\big|{\cal T}_1\big|.
\end{equation}
Similarly to (1) of the proof of Lemma~\ref{cusp2_lmm2},
\begin{gather*}\begin{split}
N(\tilde{\al})&=\blr{c\big(\ga_0^*\!\otimes\!\ga_1^{*\otimes d}\oplus
\ga_0^*\!\otimes\!\ga_1^{*\otimes d}\!\otimes\!\ga^*\big)
c(\ga^{\otimes2})^{-1},{\cal K}_1'}
-{\cal C}_{\tilde{\al}^{-1}(0)}(\tilde{\al}^{\perp})\\
&=\lr{(3\la\!+\!2y\!+\!2da)\big(\la^2\!+\!(2y\!+\!(2d\!-\!3)a)\la+
(2d\!-\!3)ya\!+\!(d^2\!-\!3d\!+\!3)a^2\big),
\Bbb{P}T\PP|_{{\cal N}_1''}}-{\cal C}_{{\cal T}_1}(\tilde{\al}^{\perp})\\
&=\lr{7y^2+(14d\!-\!39)ya+(7d^2\!-\!39d\!+\!54)a^2,{\cal N}_1''}
-{\cal C}_{{\cal T}_1}(\tilde{\al}^{\perp})\\
&=\big(7|{\cal N}_1|+(14d\!-\!39)|{\cal N}_{1,1}|+(7d^2\!-\!39d\!+\!54)\big)
-{\cal C}_{{\cal T}_1}(\tilde{\al}^{\perp}),
\end{split}\\
\hbox{where}\qquad {\cal C}_{\partial M}(\tilde{\al}^{\perp})=\big|{\cal T}_1\big|.
\end{gather*}
Combining these observations with \e_ref{cusp2gen_e1} and \e_ref{cusp2gen_e2},
we obtain
\begin{equation*}\begin{split}
|{\cal K}_2|  &=  3(d\!-\!1)^2|{\cal K}_1|
-2\big(7|{\cal N}_1|+(14d\!-\!39)|{\cal N}_{1,1}|+(7d^2\!-\!39d\!+\!54)\big)
-2|{\cal T}_1|\\
&=12(d\!-\!3)\big(3d^3\!-\!6d^2\!-\!11d\!+\!18\big).
\end{split}\end{equation*}

\subsubsection{The Number ${\cal N}_3$}
\label{n3gen_subsub}

\noindent
We take ${\cal D}\!\approx\!\Bbb{P}^3$ as above and define
${\cal N}_1''$, $\tilde{\cal N}_{2;0}'$, $\tilde{\cal N}_2'$,
$M$, $M^0$, $\tilde{\cal N}_3$, ${\cal Z}_{k,i}$ for $k\!=\!0,1$ and $i\!\!=\!1,2$,
${\cal Z}_{0,12}$, and $\al$ as in Lemmas~\ref{node3_lmm0}
and Lemma~\ref{node3_lmm1}, except now
\begin{gather*}
\vph_2\in\Ga\big({\cal N}_1''\!\times\!\PP_2;\ga_0^*\!\otimes\!\ga_2^{*\otimes d}
\oplus \ga_0^*\!\otimes\!\ga_2^{*\otimes d}\!\otimes\!T^*\PP_2\big), \qquad
\vph_3\in\Ga\big(M;\ga_0^*\!\otimes\!\ga_3^{*\otimes d}\oplus
\ga_0^*\!\otimes\!\ga_3^{*\otimes d}\!\otimes\!T^*\PP_3\big),\\
\hbox{and}\qquad
\al\in\Ga\big(\tilde{\cal N}_2';\Hom(T\PP_1,\ga_0^*\!\otimes\!\ga_1^{*\otimes d}\oplus
\ga_0^*\!\otimes\!\ga_1^{*\otimes d}\!\otimes\!T^*\PP_1)\big).
\end{gather*}
Since $\vph_3|_{M^0}$ is transverse to the zero set,
\begin{equation}\label{node3gen_e1}\begin{split}
|\tilde{\cal N}_3| &=\, ^{\pm}\big|\vph_3^{-1}(0)\!\cap\!M^0\big|
= \blr{e(\ga_0^*\!\otimes\!\ga_3^{*\otimes d}\oplus
\ga_0^*\!\otimes\!\ga_3^{*\otimes d}\!\otimes\!T^*\PP_3),M}
-{\cal C}_{\partial M}(\vph_3)\\
&= \blr{(y\!+\!da_3)(y^2\!+\!(2d-3)ya_3\!+\!(d^2\!-\!3d\!+\!3)a_3^2),\tilde{\cal N}_2'\!\times\!\PP_3}
-{\cal C}_{\partial M}(\vph_3)\\
&= 3(d\!-\!1)^2\lr{y,\tilde{\cal N}_2'}-{\cal C}_{\partial M}(\vph_3)
=6(d\!-\!1)^2|{\cal N}_2|-2{\cal C}_{{\cal Z}_{1,1}}(\vph_3)
-2{\cal C}_{{\cal Z}_{0,1}}(\vph_3)-{\cal C}_{{\cal Z}_{0,12}}(\vph_3).
\end{split}\end{equation}
Similarly to the proof of Lemma~\ref{node3_lmm1}, we have
\begin{equation}\label{node3gen_e2}
{\cal C}_{{\cal Z}_{1,1}}(\vph)=N(\al), \qquad 
{\cal C}_{{\cal Z}_{0,1}}(\vph)=3\big|{\cal K}_2\big|,
\quad\hbox{and}\quad
{\cal C}_{{\cal Z}_{0,12}}(\vph)=4\big|{\cal T}_1\big|.
\end{equation}
In order to compute $N(\al)$, we define $M$, ${\cal Z}_{0,1}$, and
${\cal Z}_{0,2}$ as in the proof of Lemma~\ref{node3_lmm2}.
By the same argument as before, we find that
\begin{gather*}\begin{split}
N(\al)&=\blr{c\big(\ga_0^*\!\otimes\!\ga_1^{*\otimes d}\oplus
\ga_0^*\!\otimes\!\ga_1^{*\otimes d}\!\otimes\!T^*\PP_1\big)c(T\PP_1)^{-1},\tilde{\cal N}_2'}
-{\cal C}_{\tilde{\al}^{-1}(0)}(\tilde{\al}^{\perp})\\
&=\lr{3y+(3d\!-\!2)a_1,\tilde{\cal N}_2'}-{\cal C}_{\partial M}(\tilde{\al}^{\perp})
=\big(6|{\cal N}_2|+3(d\!-\!2)|{\cal N}_{2,1}|\big)
-{\cal C}_{{\cal Z}_{0,1}}(\tilde{\al}^{\perp})-{\cal C}_{{\cal Z}_{0,2}}(\tilde{\al}^{\perp}),
\end{split}\\
\hbox{where}\qquad
{\cal C}_{{\cal Z}_{0,1}}(\tilde{\al}^{\perp})=|{\cal K}_2| 
\quad\hbox{and}\quad
{\cal C}_{{\cal Z}_{0,2}}(\tilde{\al}^{\perp})=2|{\cal T}_1|.
\end{gather*}
Combining this result with \e_ref{node3gen_e1} and \e_ref{node3gen_e2}, we conclude that
\begin{equation*}\begin{split}
|\tilde{\cal N}_3|&=
6\big((d\!-\!1)^2\!-\!2)|{\cal N}_2|-6(d\!-\!2)|{\cal N}_{2,1}|-4|{\cal K}_2|\\
&=3(9d^6\!-\!54d^5\!+\!9d^4\!+\!423d^3\!-\!458d^2\!-\!829d\!+\!1050).
\end{split}\end{equation*}
{\it Remark:} For $d\!\ge\!5$, the middle component of the polynomial $\al_0$ in 
the proof of Lemma~\ref{node3_lmm1} should be increased by~$\frac{1}{4}\be_5v^5$.
However, this term vanishes as we proceed from \e_ref{node3_e5} to~\e_ref{node3_e6}.

\section{Stable Maps and Recursive Formula}
\label{all_deg_sec}

\subsection{The Moduli Space of Four Marked Points on a Sphere}
\label{m04_subs}

\noindent
In this section, we derive recursion~\e_ref{recursion_e_main},
following the argument in~\cite{RuT}.
We start by defining an invariant that counts holomorphic maps into~$\P$.
A priori, the number we describe depends on the cross ratio of 
the chosen four points on a sphere.
However, it turns out that this number is well-defined.
We use its independence to express this invariant in terms
of the numbers~$n_d$ in two different ways.
By comparing the two expressions, we obtain~\e_ref{recursion_e_main}.\\

\noindent
Let $x_0,x_1,x_2$ and $x_3$ be the four points in $\PP$ given~by
$$x_0=[1,0,0],\quad x_1=[0,1,0],\quad
x_2=[0,0,1],\quad x_3=[1,1,1].$$
We denote by $H^0(\PP;\ga^{*\otimes2})$
the space of holomorphic sections of the holomorphic line bundle 
$\ga^{*\otimes2}\!\lra\!\PP$, or equivalently of
the degree-two homogeneous polynomials in three variables;
see Lemma~\ref{holom_secs_lmm}.
Let
\begin{equation*}\begin{split}
{\cal U}&=\big\{([s],x)\!\in\!
\Bbb{P}H^0(\PP;\ga^{*\otimes2})\!\times\!\PP:
s(p_i)\!=\!0 ~\forall i\!=\!0,1,2,3;~ s(x)\!=\!0\big\}\\
&\approx\big\{\big([A,B];[z_0,z_1,z_2]\big)\!\in\!
\Bbb{P}^1\!\times\!\PP\!: Az_0z_1\!-\!(A\!-\!B)z_0z_2\!-\!Bz_1z_2\!=\!0\big\}.
\end{split}\end{equation*}
The space ${\cal U}$ is a compact complex two-manifold.\\

\noindent
Let $\pi\!:{\cal U}\!\lra\!\ov{\frak M}_{0,4}\!\equiv\!\Bbb{P}^1$ 
denote the projection onto the first component.
If $[A,B]\!\in\!\ov{\frak M}_{0,4}$, the fiber~$\pi^{-1}([A,B])$ is the~conic
$${\cal C}_{A,B}=\big\{[z_0,z_1,z_2]\!\in\!\PP\!:
Az_0z_1\!-\!(A\!-\!B)z_0z_1\!-\!Bz_1z_2\!=\!0\big\}.$$
If $[A,B]\!\neq\![1,0],[0,1],[1,1]$, 
${\cal C}_{A,B}$ is a smooth complex curve
of genus zero. In other words, ${\cal C}_{A,B}$ is a sphere with
four distinct marked points by Lemma~\ref{genus_lmm}.
If $[A,B]\!=\![1,0],[0,1],[1,1]$, ${\cal C}_{A,B}$~is
a union of two lines. One of the lines contains two of
the four points $x_0,\ldots,x_3$, and the other line
passes through the remaining two points.
The two lines intersect in a single point.
Figure~\ref{m04_fig} shows the three singular fibers of 
the projection map $\pi\!:{\cal U}\!\lra\!\ov{\frak M}_{0,4}$.
The other fibers are smooth conics.
The fibers should be viewed as lying in planes
orthogonal to the horizontal line in the~figure.\\

\noindent
We conclude this subsection with a few remarks concerning
the family ${\cal U}\!\lra\!\ov{\frak M}_{0,4}$.
These remarks are irrelevant for the purposes of the next
subsection and can be omitted.\\

\noindent
If $[A,B]\!\in\!\ov{\frak M}_{0,4}\!-\!\{[1,0],[0,1],[1,1]\}$,
${\cal C}_{A,B}$ is a smooth complex curve of genus~zero,
i.e. it is a sphere holomorphically embedded in~$\PP$.
Thus, there exists a one-to-one holomorphic map
$f\!:\Bbb{P}^1\!\lra\!{\cal C}_{A,B}$.
Using Lemma~\ref{holomor_lmm},
it can be shown directly that if $[u_i,v_i]\!=\!f^{-1}(x_i)$,
$$\frac{u_0/v_0-u_2/v_2}{u_0/v_0-u_3/v_3}:
\frac{u_1/v_1-u_2/v_2}{u_1/v_1-u_3/v_3}=\frac{B}{A}.$$
The cross-ratio is the only invariant of four distinct points
on $\Bbb{P}^1$; see~\cite{A}, for example. Thus,
\begin{gather*}
\Bbb{P}^1\!-\!\{[1,0],[0,1],[1,1]\}=
{\frak M}_{0,4}\equiv\big\{(x_0,x_1,x_2,x_3)\!\in\!(\Bbb{P}^1)^4\!:
x_i\!\neq\!x_j~\hbox{if}~i\!\neq\!j\big\}\big/\sim,\\
\hbox{where}\qquad (x_0,x_1,x_2,x_3)\sim
\big(\tau(x_0),\tau(x_1),\tau(x_2),\tau(x_3)\big)
\quad\hbox{if}\quad
\tau\!\in\!\hbox{PSL}_2\equiv\hbox{Aut}(\Bbb{P}^1).
\end{gather*}
Furthermore, the restriction of the projection map 
$\pi\!:{\cal U}|_{{\frak M}_{0,4}}\!\lra\!{\frak M}_{0,4}$ to each fiber
${\cal C}_{[A,B]}$ 
is the cross ratio of the points $x_0,\ldots,x_3$ on~${\cal C}_{[A,B]}$,
viewed as an element of $\Bbb{P}^1\!\supset\!\Bbb{C}$.

\begin{figure}
\begin{pspicture} (-1.1,-1.5)(10,3)
\psset{unit=.4cm}
\psline[linewidth=.03](2,-1)(33,-1)
\rput(30,-2){$\ov{\frak M}_{0,4}\!\approx\!\Bbb{P}^1$}
\rput(30,4){${\cal U}$}\rput(30.5,2){$\pi$}
\psline[linewidth=.03]{->}(30,3)(30,0)
\pscircle*(5,-1){.2}\pscircle*(15,-1){.2}\pscircle*(25,-1){.2}
\rput(5,-2){$[1,0]$}\rput(15,-2){$[1,1]$}\rput(25,-2){$[0,1]$}
\psline[linewidth=.02](3,6)(7,0)\psline[linewidth=.03](3,0)(7,6)
\pscircle*(5,3){.14}\pscircle*(4.5,2.25){.17}\pscircle*(5.5,3.75){.17}
\pscircle*(5.25,2.625){.17}\pscircle*(4.333,4){.17}
\rput(3.85,2.35){$x_0$}\rput(3.9,3.7){$x_1$}
\rput(6.2,3.6){$x_3$}\rput(6,2.7){$x_2$}
\psline[linewidth=.02](13,1.5)(18,4)
\psline[linewidth=.02](13,4.29)(18,3.17)
\pscircle*(16.9,3.45){.14}\pscircle*(14.5,2.25){.17}\pscircle*(15.5,3.75){.17}
\pscircle*(15.25,2.625){.17}\pscircle*(14.333,4){.17}
\rput(13.8,2.4){$x_0$}\rput(13.85,3.65){$x_1$}
\rput(16.2,3.9){$x_3$}\rput(15.9,2.4){$x_2$}
\psline[linewidth=.02](24.233,5.05)(24.75,-.275)
\psline[linewidth=.02](25.7,4.65)(24.55,-.525)
\pscircle*(24.7,.15){.14}\pscircle*(24.5,2.25){.17}\pscircle*(25.5,3.75){.17}
\pscircle*(25.25,2.625){.17}\pscircle*(24.333,4){.17}
\rput(23.8,2.4){$x_0$}\rput(23.85,3.65){$x_1$}
\rput(26.2,3.9){$x_3$}\rput(26,2.6){$x_2$}
\end{pspicture}
\caption{The Family ${\cal U}\!\lra\!\ov{\frak M}_{0,4}$}
\label{m04_fig}
\end{figure}

\subsection{Counts of Holomorphic Maps}
\label{gw_subs}

\noindent
If $d$ is an integer and ${\cal C}$ is a complex curve, 
which may be a wedge of spheres, let
\begin{equation}\label{recursion_e0}
{\cal H}_d({\cal C})= \big\{ f\!\in\!C^{\i}({\cal C};\PP)\!:
f~\hbox{\it is holomorphic},~f_*[{\cal C}]\!=\!d[L]\big\},
\end{equation}
where $[L]\!\in\!H_2(\PP;\Bbb{Z})$ is the homology class
of a line in~$\PP$.
We give a more explicit description of the space ${\cal H}_d({\cal C})$
in the relevant cases~below.\\

\noindent
Suppose $\ell_0,\ell_1$ and $p_2,\ldots,p_{3d-1}$ are 
two lines and $3d\!-\!2$ points in general position in~$\PP$. 
If $\si\!\in\!\ov{\frak M}_{0,4}$, let 
$N_d^{\si}(l_0,l_1,p_2,\ldots,p_{3d-1})$ denote the cardinality of the~set
\begin{equation}\label{recursion_e1}
\big\{f\!\in\!{\cal H}_d({\cal C}_{\si})\!:
f(x_0)\!\in\!\ell_0,~f(x_1)\!\in\!\ell_1,~
f(x_2)\!=\!p_2,~f(x_3)\!=\!p_3,~p_i\!\in\!\hbox{Im}~\!f~\forall i\big\}.
\end{equation}
Here ${\cal C}_{\si}$ denotes the rational curve with four marked points,
$x_0$, $x_1$, $x_2$, and $x_3$, 
whose cross ratio is~$\si$; see Subsection~\ref{m04_subs}.
If $\si\!\neq\![1,0],[0,1],[1,1]$, 
${\cal C}_{\si}$~is a sphere with four, distinct, marked points.
In this case, the condition $f\!\in\!{\cal H}_d({\cal C}_{\si})$ means 
that $f$ has the form
$$f([u,v])=\big[P_0(u,v),P_1(u,v),P_2(u,v)\big]\qquad
\forall[u,v]\!\in\!\Bbb{P}^1,$$
for some degree-$d$ homogeneous polynomials $P_0,P_1,P_2$ that
have no common factor; see Lemma~\ref{holomor_lmm}.
If $\si\!=\![1,0],[0,1],[1,1]$, ${\cal C}_{\si}$ is a wedge of two spheres, 
${\cal C}_{\si,1}$ and~${\cal C}_{\si,2}$, with two marked points each.
In this case, the first condition in~\e_ref{recursion_e0} means
that $f$ is continuous and $f|_{{\cal C}_{\si,1}}$ and
$f|_{{\cal C}_{\si,2}}$ are holomorphic.
The second condition in~\e_ref{recursion_e0} means that $d\!=\!d_1\!+\!d_2$ if
$f_*[{\cal C}_{\si,1}]\!=\!d_1[L]$ and
$f_*[{\cal C}_{\si,2}]\!=\!d_2[L]$.\\

\noindent
The requirement that the two lines, $\ell_0$ and $\ell_1$,
and the $3d\!-\!2$ points, $p_2,\ldots,p_{3d-1}$, are in general position
means that they lie in a dense open subset~${\cal U}_{\si}$ of the space of all
possible tuples $(\ell_0,\ell_1,p_2,\ldots,p_{3d-1})$:
$${\frak X}\equiv
\hbox{Gr}_2\Bbb{C}^3\times \hbox{Gr}_2\Bbb{C}^3\times\big(\PP\big)^{3d-2}.$$
Here $\hbox{Gr}_2\Bbb{C}^3$ denotes the Grassmanian manifold of
two-planes through the origin in~$\Bbb{C}^3$, or equivalently of
lines in~$\PP$.
The dense open subset ${\cal U}_{\si}$ of ${\frak X}$ consists of
tuples $(\ell_0,\ell_1,p_2,\ldots,p_{3d-1})$
that satisfy a number of geometric conditions.
In particular, $\ell_0\!\neq\!\ell_1$, none of the points $p_2,\ldots,p_{3d-1}$
lies on either $\ell_0$ or $\ell_1$, 
the $3d\!-\!1$ points $\ell_0\cap\ell_1, p_2,\ldots,p_{3d-1}$
are distinct, no three of them lie on the same line, and so on.
In addition, we need to impose certain cross-ratio conditions
on the rational curves that pass through $\ell_0$, $\ell_1$,
$p_2$, $p_3$, and a subset of the remaining $3d\!-\!4$ points.
These conditions can be stated more formally. Define
$$\hbox{ev}_{\si}\!:
{\cal H}_d({\cal C}_{\si})\times\big({\cal C}_{\si}\big)^{3d-4}\lra
\big(\PP)^d\qquad\hbox{by}\qquad
\hbox{ev}_{\si}\big(f;x_4,\ldots,x_{3d-1}\big)=
\big(f(x_1),\ldots,f(x_{3d-1})\big).$$
Lemma~\ref{holomor_lmm} implies that ${\cal H}_d({\cal C}_{\si})$ 
is a dense open subset of~$\Bbb{P}^{3d+2}$ and
the evaluation map~$\hbox{ev}_{\si}$ is holomorphic.
The space ${\cal H}_d({\cal C}_{\si})$ has a natural 
compactification~$\ov{\frak M}_{\si}(\PP,d)$,
which is the union of spaces of holomorphic maps from 
various wedges of spheres into~$\PP$.
The complex dimension of each such boundary stratum is less
than that of~${\cal H}_d({\cal C}_{\si})$.
The evaluation map~$\hbox{ev}_{\si}$ admits a continuous extension over
$\partial\ov{\frak M}_{\si}(\PP,d)$, whose restriction to
each stratum is holomorphic.
The elements $(\ell_0,\ell_1,p_2,\ldots,p_{3d-1})$
of the subspace~${\cal U}_{\si}$ of~${\frak X}$ 
are characterized by the condition that 
the restriction of the evaluation map to each stratum of
$\ov{\frak M}_{\si}(\PP,d)$ is transversal to the submanifold
$$\ell_0\times\ell_1\times\!p_2\!\times\!\ldots\!\times\!p_{3d-1}
\subset (\PP)^{3d}.$$
This condition implies that
$$\hbox{ev}_{\si}^{-1}\big(\ell_0\times\ell_1\times\!p_2\!\times\!\ldots\!\times\!p_{3d-1}\big)
\cap\partial\ov{\frak M}_{\si}(\PP,d)=\eset$$
and the set in \e_ref{recursion_e1} is a finite subset of~${\cal H}_d({\cal C}_{\si})$.\\

\noindent
The set ${\cal U}_{\si}$ of "general" tuples $(\ell_0,\ell_1,p_2,\ldots,p_{3d-1})$
is path-connected. Indeed, it is the complement of a finite number 
of proper complex submanifolds in~${\frak X}$.
It follows that the number in~\e_ref{recursion_e1} is independent of the choice of 
two lines and $3d\!-\!2$ points in general position in~$\PP$.
We thus may simply denote it by~$N_d^{\si}$.
If $\si\!\neq\![1,0],[0,1],[1,1]$, ${\cal C}_{\si}$ is a sphere with
four distinct points.
In such a case, it is fairly easy to show that the number~$N_d^{\si}$ 
does not change with small variations~$\si$, or equivalently of
the four points on the sphere.
Thus, $N_d^{\si}$ is independent of~
$$\si\in{\frak M}_{0,4}=\Bbb{P}^1-\big\{[1,0],[0,1],[1,1]\big\}
=\ov{\frak M}_{0,4}-\big\{[1,0],[0,1],[1,1]\big\}.$$
It is far harder to prove

\begin{prp}
\label{gluing_prp}
The function $\si\!\lra\!N_d^{\si}$ is constant on $\ov{\frak M}_{0,4}$.
\end{prp}

\noindent
This proposition is a special case of the gluing theorems first
proved in~\cite{McSa} and~\cite{RuT}.
A more straightforward proof can be obtained via the approach of~\cite{LT}.

\subsection{Holomorphic Maps vs.~Complex Curves}
\label{recursion_subs}

\noindent
In this subsection, we express the numbers $N_d^{[1,0]}$ and $N_d^{[0,1]}$
of Subsections~\ref{gw_subs} in terms of the numbers~$n_{d'}$,
with $d'\!\le\!d$, of Question~\ref{g0n2_ques}.
By Proposition~\ref{gluing_prp}, $N_d^{[1,0]}\!=\!N_d^{[0,1]}$.
We obtain a recursion for the numbers of~Question~\ref{g0n2_ques}
by comparing the expressions for~$N_d^{[1,0]}$ and~$N_d^{[0,1]}$.\\

\noindent
Let ${\cal C}_1$ denote the component of ${\cal C}_{[1,0]}$
containing the marked points~$x_0$ and~$x_3$;
see Figure~\ref{m04_fig}.
We denote by ${\cal C}_2$ the other component of ${\cal C}_{[1,0]}$.
By definition, 
\begin{gather*}
N_d^{[1,0]}\!=\!\sum_{d_1+d_2=d}\!N_{d_1,d_2}^{[1,0]}  \qquad\hbox{where}\\ 
\begin{split}
N_{d_1,d_2}^{[1,0]}=
\big|\big\{f\!\in\!{\cal H}_d({\cal C}_{[0,1]};\PP)\!:~&
f_*[{\cal C}_1]\!=\!d_1[L]~,f_*[{\cal C}_2]\!=\!d_2[L];~
p_i\!\in\!\hbox{Im}~\!f~\forall i;\\
& f(x_0)\!\in\!\ell_0,~f(x_1)\!\in\!\ell_1,~
f(x_2)\!=\!p_2,~f(x_3)\!=\!p_3\big\}\big|.
\end{split}
\end{gather*}
Since the group $PSL_2$ of holomorphic automorphisms acts transitively
on triples of distinct points on the sphere,
\begin{equation*}\begin{split}
N_{d_1,d_2}^{[1,0]}=
\big|\big\{(f_1,f_2)\!\in\!{\cal H}_{d_1}(S^2)\!\times\!{\cal H}_{d_2}(S^2)\!:~&
f_1(\i)\!=\!f_2(\i),~p_i\!\in\!f_1(S^2)\cup f_2(S^2)~\forall i;\\
&f_1(0)\!\in\!\ell_0,~f_1(1)\!=\!p_3,~f_2(0)\!\in\!\ell_1,~f_2(1)\!=\!p_2\big\}\big|.
\end{split}\end{equation*}
Since the maps $f_1$ and $f_2$ above are holomorphic, 
$d_1,d_2\!\ge\!0$ if $N_{d_1,d_2}^{[1,0]}\!\neq\!0$.
Since every degree-zero holomorphic map is constant and $p_3\!\not\in\!\ell_0$,
$N_{0,d}^{[1,0]}\!=\!0$.
Similarly, $N_{d,0}^{[1,0]}\!=\!0$. Thus, we assume that $d_1,d_2\!>\!0$.
Since the points $p_3,\ldots,p_{3d-1}$ are in general position, 
$f_1(S^2)$ contains at most $3d_1\!-\!2$ of the points $p_4,\ldots,p_{3d-1}$.
Similarly, the curve $f_2(S^2)$ passes through at most $3d_2\!-\!2$
of the points $p_4,\ldots,p_{3d-1}$. Thus, if $I\!=\!\{4,\ldots,3d\!-\!1\}$,
$$N_{d_1,d_2}^{[1,0]}=\!\!\!\sum_{I=I_1\sqcup I_2,|I_1|=3d_1-2}\!\!\!\!\!\!
N_{d_1,d_2}^{[1,0]}(I_1,I_2),$$
where $N_{d_1,d_2}^{[1,0]}(I_1,I_2)$ is the cardinality of the~set
\begin{equation*}\begin{split}
{\cal S}_{d_1,d_2}^{[1,0]}(I_1,I_2)=
\big\{(f_1,f_2)\!\in\!{\cal H}_{d_1}(S^2)\!\times\!{\cal H}_{d_2}(S^2)\!:~
p_i\!\in\!f_1(S^2)~\forall i\!\in\!I_1,~
p_i\!\in\!f_2(S^2)~\forall i\!\in\!I_2;&\\
f_1(\i)\!=\!f_2(\i),~f_1(0)\!\in\!\ell_0,~f_1(1)\!=\!p_3,
~f_2(0)\!\in\!\ell_1,~f_2(1)\!=\!p_2&\big\}.
\end{split}\end{equation*}
If $(f_1,f_2)\!\in\!{\cal S}_{d_1,d_2}^{[1,0]}(I_1,I_2)$, $f_1(S^2)$
is one of the $n_{d_1}$ curves passing through the points
$\{p_i\!:i\!\in\!\{3\}\sqcup I_1\}$.
Similarly, $f_2(S^2)$ is one of the $n_{d_2}$ curves passing through 
the points $\{p_i\!:i\!\in\!\{2\}\sqcup I_2\}$.
The point $f_1(\i)\!=\!f_2(\i)$ must be one of the $d_1d_2$ 
points of $f_1(S^2)\cap f_2(S^2)$; see Lemma~\ref{bezout_lmm}.
Finally, $f_1(0)$ must be one of the $d_1$~points of $f_1(S^2)\cap\ell_0$,
while $f_2(0)$ must be one of the $d_2$~points of $f_2(S^2)\cap\ell_1$.
Thus, we conclude that
\begin{equation}\label{recursion_e5}\begin{split}
N_d^{[1,0]}=\!\!\sum_{d_1+d_2=d}\!\!\!N_{d_1,d_2}^{[1,0]}
&=\!\!\!\sum_{d_1+d_2=d}\sum_{I=I_1\sqcup I_2,|I_1|=3d_1-2}\!\!\!\!\!\!
N_{d_1,d_2}^{[1,0]}(I_1,I_2)\\
&=\!\!\sum_{d_1+d_2=d}\sum_{I_1\subset I,|I_1|=3d_1-2}\!\!\!\!\!\!\!\!\!
(d_1d_2)(d_1n_{d_1})(d_2n_{d_2})\\
&=\sum_{d_1+d_2=d}\binom{3d\!-\!4}{3d_1\!-\!2}d_1^2d_2^2n_{d_1}n_{d_2};
\end{split}\end{equation}
where $I\!=\!\{4,\ldots,3d\!-\!1\}$.\\

\noindent
We compute the number $N_d^{[0,1]}$ similarly.
We denote by ${\cal C}_1$ the component of ${\cal C}_{[0,1]}$
containing the points $x_0$ and~$x_1$ and 
by ${\cal C}_2$ the other component of ${\cal C}_{[0,1]}$.
By definition,
\begin{gather*}
N_d^{[0,1]}\!=\!\sum_{d_1+d_2=d}\!N_{d_1,d_2}^{[0,1]},  \qquad\hbox{where}\\ 
\begin{split}
N_{d_1,d_2}^{[0,1]}=
\big|\big\{(f_1,f_2)\!\in\!{\cal H}_{d_1}(S^2)\!\times\!{\cal H}_{d_2}(S^2)\!:~&
f_1(\i)\!=\!f_2(\i),~p_i\!\in\!f_1(S^2)\cup f_2(S^2)~\forall i;\\
&f_1(0)\!\in\!\ell_0,~f_1(1)\!\in\!\ell_1,~f_2(0)\!=\!p_2,~f_2(1)\!=\!p_3\big\}\big|.
\end{split}
\end{gather*}
Since every degree-zero holomorphic map is constant, $N^{[0,1]}_{d,0}\!=\!0$
as before. However, 
$$N^{[0,1]}_{0,d}=\big|\big\{f_2\!\in\!{\cal H}_d(S^2)\!:
f_2(\i)\!\in\!\ell_0\cap\ell_1,~p_i\!\in\!f_2(S^2)~\forall i\!=\!2.\ldots,3d\!-\!1\big\}\big|
=n_d.$$
If $d_1,d_2\!>\!0$, 
$$N_{d_1,d_2}^{[0,1]}=\!\!\!\sum_{I=I_1\sqcup I_2,|I_1|=3d_1-1}\!\!\!\!\!\!
N_{d_1,d_2}^{[0,1]}(I_1,I_2),$$
where $N_{d_1,d_2}^{[0,1]}(I_1,I_2)$ is the cardinality of the~set
\begin{equation*}\begin{split}
{\cal S}_{d_1,d_2}^{[0,1]}(I_1,I_2)=
\big\{(f_1,f_2)\!\in\!{\cal H}_{d_1}(S^2)\!\times\!{\cal H}_{d_2}(S^2)\!:~
p_i\!\in\!f_1(S^2)~\forall i\!\in\!I_1,~
p_i\!\in\!f_2(S^2)~\forall i\!\in\!I_2;&\\
f_1(\i)\!=\!f_2(\i),~f_1(0)\!\in\!\ell_0,~f_1(1)\!\in\!\ell_1,
~f_2(0)\!=\!p_2,~f_2(1)\!=\!p_3&\big\}.
\end{split}\end{equation*}
Proceeding as in the previous paragraph, we conclude that
\begin{equation}\label{recursion_e7}\begin{split}
N_d^{[0,1]}=\!\!\sum_{d_1+d_2=d}\!\!\!N_{d_1,d_2}^{[0,1]}
&=n_d+\!\sum_{d_1+d_2=d}\sum_{I=I_1\sqcup I_2,|I_1|=3d_1-1}\!\!\!\!\!\!
N_{d_1,d_2}^{[0,1]}(I_1,I_2)\\
&=n_d+\!\sum_{d_1+d_2=d}\sum_{I_1\subset I,|I_1|=3d_1-1}\!\!\!\!\!\!\!\!\!
(d_1d_2)(d_1^2n_{d_1})(n_{d_2})\\
&=n_d+\sum_{d_1+d_2=d}\binom{3d\!-\!4}{3d_1\!-\!1}d_1^3d_2n_{d_1}n_{d_2};
\end{split}\end{equation}
Comparing equations~\e_ref{recursion_e5} and~\e_ref{recursion_e7}, 
we obtain
\begin{equation}\label{recursion_e9}
n_d=\sum_{d_1+d_2=d}\Bigg(
\binom{3d\!-\!4}{3d_1\!-\!2}d_1d_2-\binom{3d\!-\!4}{3d_1\!-\!1}d_1^2\Bigg)
d_1d_2n_{d_1}n_{d_2}.
\end{equation}
The recursive formula~\e_ref{recursion_e_main} is the symmetrized
version of~\e_ref{recursion_e9}.

\appendix

\section{The Basics}
\label{basics_sec}

\subsection{Complex Projective Spaces}
\label{projective_subs}

\noindent
The {\it complex projective space} $\P$ is the space of (complex) lines
through the origin in~$\Bbb{C}^{n+1}$. Equivalently,
$$\P=\big(C^{n+1}\!-\!\{0\}\big)\big/C^*, \qquad \hbox{where}\quad
\big(z_0,\ldots,z_n\big)\sim\big(tz_1,\ldots,tz_n)
~~\hbox{if}~t\!\in\!\Bbb{C}^*.$$
This space is a smooth $2n$-manifold. For $i\!=\!0,\ldots,n$, let
\begin{gather*}
U_i=\big\{[z_0,\ldots,z_n]\!\in\!\P\!: z_i\!\neq\!0\big\},\\
\phi_i\!:\Bbb{C}^n\lra U_i,\quad
\phi_i(w_1,\ldots,w_n)=\big[w_1,\ldots,w_i,1,w_{i+1},\ldots,w_n\big].
\end{gather*}
The set $\big\{(U_i,\phi_i,\Bbb{C}^n)\big\}$ is the {\it standard atlas}
for~$\P$. If $i\!<\!j$, the corresponding overlap map is given~by
\begin{gather*}
\phi_{ij}\!\equiv\!\phi_i^{-1}\circ\phi_j\big|_{\phi_j^{-1}(U_i)}\!:
\big\{(w_1,\ldots,w_n)\!\in\!\Bbb{C}^n\!: w_{i+1}\!\neq\!0\big\}\lra
\big\{(w_1,\ldots,w_n)\!\in\!\Bbb{C}^n\!: w_j\!\neq\!0\big\}\\
(w_1,\ldots,w_n)\lra
\Big(\frac{w_1}{w_{i+1}},\ldots,\frac{w_i}{w_{i+1}},\frac{w_{i+2}}{w_{i+1}},\ldots,
\frac{w_j}{w_{i+1}},w_{i+1}^{-1},\frac{w_{j+1}}{w_{i+1}},\ldots,\frac{w_n}{w_{i+1}}\Big).
\end{gather*}
Each map~$\phi_{ij}$ is a diffeomorphism.
In fact, this map is holomorphic, and so is its inverse~$\phi_{ij}^{-1}$.
In other words, $\P$~is naturally a {\it complex $n$-manifold}.\\

\noindent
Suppose $X$ and $Y$ are complex manifolds, of complex dimensions $m$ and~$n$,
and with (holomorphic) atlases $\big\{(U_i,\phi_i,U_i')\big\}_{i\in I}$ and $\big\{(V_j,\varphi_j,V_j')\big\}_{j\in J}$, respectively.
Smooth map \hbox{$f\!:X\!\lra\!Y$} is called {\it holomorphic} if
for all $i\!\in\!I$ and $j\!\in\!J$, the~map
$$\varphi_j^{-1}\circ f\circ\phi_i\!:
\phi_i^{-1}(f^{-1}(V_j))\lra\Bbb{C}^n$$
is holomorphic as a $\Bbb{C}^n$-valued function on an open subset of~$\Bbb{C}^m$.
In the case of interest to~us, i.e.~$X\!=\!\Bbb{P}^1$ and $Y\!=\!\P$,
the holomorphic maps have a much simpler description,
see Lemma~\ref{holomor_lmm} below.
This lemma can be checked directly.
The simpler characterization of Lemma~\ref{holomor_lmm} can be taken 
as the definition of what it means to be a holomorphic map between $\Bbb{P}^1$ and~$\P$.

\begin{lmm}
\label{holomor_lmm}
If $f\!:\Bbb{P}^1\!\lra\!\P$ is a holomorphic map, there exist homogeneous
polynomials $p_0,\ldots,p_n$ in two variables
such that $p_0,\ldots,p_n$ are of the same degree,
have no common factor, and 
\begin{equation}\label{holomor_lmm_e}
f\big([z_0,z_1]\big)=\big[p_0(z_0,z_1),\ldots,p_n(z_0,z_1)\big]
\qquad\forall [z_0,z_1]\!\in\!\Bbb{P}^1.
\end{equation}
Conversely, if $p_0,\ldots,p_n$ are homogeneous polynomials in two variables 
that are of the same degree and have no common factor, the map $f\!:\Bbb{P}^1\!\lra\!\P$ given 
by~\e_ref{holomor_lmm_e} is well-defined and holomorphic.
\end{lmm}

\subsection{Almost Complex and Symplectic Structures}
\label{sympl_subs}

\noindent
This subsection is not relevant for understanding Sections~\ref{low_deg_sec}-\ref{all_deg_sec}. 
However, it puts the last section in perspective.\\

\noindent
Let $X$ be a smooth manifold. An {\it almost complex structure} on~$X$
is a smooth section~$J$ of the bundle $\hbox{End}(TX)\!\lra\!X$ 
such that $J^2\!=\!-I$.
In other words, an almost complex structure is a smooth family of
linear maps $J_p\!\!:T_pX\!\lra\!T_pX$ such that $J_pJ_pv\!=\!-v$
for all $v\!\in\!T_pX$ and $p\!\in\!X$.
For example, if $X\!=\!\Bbb{C}^n$, $T_p\Bbb{C}^n\!=\!\Bbb{C}^n$ 
and the desired endomorphism on $T_p\Bbb{C}^n$ is simply 
the multiplication by~$i$.\\

\noindent
Every complex $n$-manifold~$X$ carries a natural almost complex structure~$J$,
defined as follows.
Let $\big\{(U_i,\phi_i,U_i')\big\}_{i\in I}$ be the (holomorphic) atlas for~$X$.
If $p\!\in\!U_i$, we~set
$$J_p=d\phi_i\big|_{\phi_i^{-1}(p)}\circ i\circ d\phi_i^{-1}\big|_p.$$ 
Since all overlap maps $\phi_i^{-1}\!\circ\!\phi_j$ are holomorphic,
the endomorphism~$J_p$ is independent of the choice of $i\!\in\!I$
such that $p\!\in\!U_i$.
An almost complex structure arising in such a way is called {\it complex}
or {\it integrable}.\\

\noindent
A typical almost complex structure is not integrable, 
unless the real dimension of the manifold is~two.
In fact, there is a criterion that characterizes integrable 
almost complex structures.
If $(X,J)$ is an almost complex manifold, $p\!\in\!X$,
and $V$ and $W$ are vector fields on~$X$, let
$$N_p^J(V_p,W_p)=\frac{1}{4}
\big([V,W]_p+J_p[JV,W]_p+J_p[V,JW]_p-[JV,JW]_p\big).$$
The vector $N_p^J(V_p,W_p)\!\in\!T_pX$ depends only on 
the values $V_p$ and $W_p$ of the vector fields $V$ and~$W$
at the point~$p$.
In addition, $N_p^J$~is linear in each of the two inputs.
Thus,
$$N^J\in\Ga\big(X;\Hom(TX\!\otimes\!TX,TX)\big),$$
i.e.~$N^J$ is a $(2,1)$-tensor field on~$X$.
This tensor field is called the {\it Nijenhuis torsion} of~$J$.
It is easy to see that $N^J\!\equiv\!0$ if $J$~is an integrable
almost complex structure. The converse is proved in~\cite{NeNi}.
Since $N^J\!\equiv\!0$ if $(X,J)$ is an almost complex manifold
of real dimension~two, it follows every almost complex structure
on a smooth two-manifold is integrable.
Such a manifold is called a {\it Riemann surface}.\\

\noindent
Suppose $(X,j)$ and $(Y,J)$ are almost complex manifolds
and $f\!:X\!\lra\!Y$ is a smooth map.
If $z\!\in\!X$, we set
$$\bar{\partial}_{J,j}f\big|_z=df\big|_z+ J_{f(z)}\circ df\big|_z\circ j_z 
\in\Hom(T_zX,T_{f(z)}Y).$$
Note that $\bar{\partial}_{J,j}f|_z\!\circ\!j_z\!=\!
-\!J_{f(z)}\!\circ\!\bar{\partial}_{J,j}f|_z$, i.e.
the linear map $\bar{\partial}_{J,j}f\big|_z$ is
$(J,j)$-antilinear. Thus,
$$\bar{\partial}_{J,j}f\in\Ga
\big(X,\La_{J,j}^{0,1}T^*X\!\otimes\!f^*TY)\big),$$
where $\La_{J,j}^{0,1}T^*X\!\otimes\!f^*TY\!\lra\!X$ is the bundle
of $(f^*J,j)$-antilinear homomorphisms from $(TX,j)$ to~$f^*(TY,J)$.
The smooth map $f\!:X\!\lra\!Y$ is called {\it $(J,j)$-holomorphic},
or {\it pseudoholomorphic}, if $\bar{\partial}_{J,j}f\!\equiv\!0$.
If $(X,j)$ and $(Y,J)$ are complex manifolds, this definition agrees
with the one given in the previous subsection.
More generally, if $(X,j)$ is a wedge of finitely many almost complex 
manifolds $(X_l,j_l)$, we will call a continuous map $f\!:X\!\lra\!Y$
{\it $(J,j)$-holomorphic} if $f|_{X_l}$ is $(J,j_l)$-holomorphic for all~$l$.\\

\noindent
If $(X,J)$ is an almost complex manifold, $A\!\in\!H_2(X;\Bbb{Z})$,
and $g$ and $n$ are nonnegative integers, let
\begin{equation*}\begin{split}
{\frak M}_{g,n}(X,A;J)=
\big\{(\Si,j,x_1,\ldots,x_n;f)\!: (\Si,j)=\hbox{Riemann surface of genus}~g;\quad&\\
x_i\!\in\!\Si,~x_i\!\neq\!x_j~\hbox{if}~i\!\neq\!j;~
f\!\in\!C^{\i}(\Si;X),~f_*[\Si]\!=\!A,~\bar{\partial}_{J,j}f\!=\!0&\big\}\big/~,
\end{split}\end{equation*}
$$\hbox{where}\quad
(\Si,j,z_1,\ldots,z_n;f)\sim \big(\Si',j',\tau(z_1),\ldots,\tau(z_n),f\!\circ\!\tau^{-1}\big)
\quad\hbox{if}\quad
\tau\!\in\!C^{\i}(\Si;\Si'),~\bar{\partial}_{j,j'}\tau\!=\!0.$$
This moduli space has a natural topology, as well as $n$ evaluation maps
$$\hbox{ev}_i\!:{\frak M}_{g,n}(X,A;J)\lra X,\quad
[\Si,j,z_1,\ldots,z_n;f]\lra f(z_i).$$
In general, ${\frak M}_{g,n}(X,A;J)$ is not a compact topological space.
However, under certain conditions on~$(X,J)$, ${\frak M}_{g,n}(X,A;J)$
admits a natural compactification and in fact carries a (virtual) 
fundamental class.\\

\noindent
Let $X$ be a smooth manifold. A {\it symplectic form} on $X$ is a closed 
two-form $\om$ on $X$ which is nondegenerate at every point of~$X$.
In other words, $d\om\!=\!0$, and for every point~$p$ in~$X$
and nonzero tangent vector $v\!\in\!T_pX$, there exists $w\!\in\!T_pX$
such that $\om_p(v,w)\!\neq\!0$.
For example, if $(x_1,y_1,\ldots,x_n,y_n)$ are the standard coordinates
on~$\Bbb{C}^n$,
$$\om\equiv dx_1\wedge dy_1+\ldots+dx_n\wedge dy_n$$
is a symplectic form on $\Bbb{C}^n$.
More generally, if $X$ admits a symplectic form, 
the (real) dimension of~$X$ is even.\\

\noindent
If $(X,\om)$ is a symplectic manifold, the almost complex structure~$J$
on~$X$ is {\it $\om$-tame} if for every point~$p$ in~$X$
and nonzero tangent vector $v\!\in\!T_pX$, $\om_p(v,J_pv)\!>\!0$.
The $\om$-tame almost complex structure~$J$ is {\it $\om$-compatible} if 
$$\om_p(Jv,Jw)\!=\!\om_p(v,w)\qquad\forall~p\!\in\!X,~v,w\!\in\!T_pX.$$ 
For example, if $\om$ is the standard symplectic form on~$\Bbb{C}^n$,
defined in the previous paragraph,
the standard complex structure~$i$, defined in the second paragraph
of this subsection, is $\om$-compatible.
For a general symplectic manifold $(X,\om)$, the spaces of $\om$-tame
and $\om$-compatible almost complex structures on~$X$ are non-empty
and contractible.
The most fundamental result in the theory of pseudoholomorphic curves
is Gromov's Compactness Theorem, stated roughly below.

\begin{thm}{\cite{Gro}}
\label{gromov_thm}
Suppose $(X,\om)$ is a compact symplectic manifold and $J$ is an almost complex
$\om$-tame structure on~$X$.
If $A\!\in\!H_2(X;\Bbb{Z})$ and $g$ and $n$ are nonnegative integers,
the moduli space ${\frak M}_{g,n}(X,A;J)$ admits a natural compactification
$\ov{\frak M}_{g,n}(X,A;J)$.
In particular, the evaluation maps~$\hbox{ev}_i$ extend continuously 
over $\ov{\frak M}_{g,n}(X,A;J)$.
\end{thm}

\noindent
The compactification $\ov{\frak M}_{g,n}(X,A;J)$ consists of equivalence classes
of tuples $(\Si,j,x_1,\ldots,x_n,f)$, where $(\Si,j)$ is a possibly singular 
genus-$g$ Riemann surface, i.e.~a wedge of smooth Riemann surfaces,
$x_1,\ldots,x_n$ are distinct points on~$\Si$, and
$f\!:\Si\!\lra\!X$ is a $(J,j)$-holomorphic map such that $f_*[\Si]\!=\!A$.
Notice that the space $\ov{\frak M}_{g,n}(X,A;J)$ is described by  
the almost complex structure~$J$, and not the symplectic form~$\om$.
However, this space may not be compact if $J$ is not $\om$-tame
for some symplectic form~$\om$ on~$X$.\\

\noindent
Since the space of $\om$-tame almost complex structures on $X$ is contractible,
up to an appropriate equivalence, the space $\ov{\frak M}_{g,n}(X,A;J)$ 
is independent of the choice of~$J$.
In particular, the "equivalence class" of $\ov{\frak M}_{g,n}(X,A;J)$
is determined by $(X,\om)$ and thus is a symplectic invariant.
This is essentially the {\it Gromov-Witten invariant} of~$(X,\om)$.

\subsection{Tautological Line Bundle}
\label{taut_subs}

\noindent
We continue with the notation of Subsection~\ref{projective_subs}.
Let
$$\ga=\big\{(\ell;z_0,\ldots,z_n)\!\in\!\P\!\times\!\Bbb{C}^{n+1}\!:
(z_0,\ldots,z_n)\!\in\!\ell\big\}.$$
We denote by $\pi\!:\ga\!\lra\!\P$ the projection map. 
For each $\ell\!\in\!\P$, the fiber $\ga_{\ell}\!\equiv\!\pi^{-1}(\ell)$
over a point $\ell\!\in\!\P$ is the line~$\ell$ through the origin
in~$\Bbb{C}^n$.
For each $i\!=\!0,\ldots,n$, let
\begin{gather*}
\tilde{U}_i=\pi^{-1}(U_i)=\big\{(\ell;z_0,\ldots,z_n)\!\in\!\ga\!:
z_i\!\neq\!0\big\},\\
\tilde{\phi}_i\!:\Bbb{C}^n\!\times\!\Bbb{C}\lra\tilde{U}_i,\quad
\tilde{\phi}_i\big(w_1,\ldots,w_n;\la)=
\big(\phi_i(w_1,\ldots,w_n);
\la w_1,\ldots,\la w_i,\la,\la w_{i+1},\ldots,\la w_n\big).
\end{gather*}
The set $\big\{(\tilde{U}_i,\tilde{\phi}_i,\Bbb{C}^n)\!\times\!\Bbb{C}\big\}$ is 
the {\it standard atlas} for~$\ga$. 
If $i\!<\!j$, the corresponding overlap map is given~by
$$\tilde{\phi}_{ij}\!\equiv\!
\tilde{\phi}_i^{-1}\circ\tilde{\phi}_j\big|_{\tilde{\phi}_j^{-1}(\tilde{U}_i)}\!\!:
\phi_j^{-1}(U_i)\!\times\!\Bbb{C}\lra\phi_i^{-1}(U_j)\!\times\!\Bbb{C},~~
(w_1,\ldots,w_n;\la)\lra\big(\phi_{ij}(w_1,\ldots,w_n);w_{i+1}\la\big).$$
Each map~$\tilde{\phi}_{ij}$ is holomorphic, and so is its inverse~$\tilde{\phi}_{ij}^{-1}$.
Thus, $\ga$~is a complex $(n\!+\!1)$-manifold.
Furthermore, if $p\!:\Bbb{C}^n\!\times\!\Bbb{C}\!\lra\!\Bbb{C}^n$ is the projection~map,
$$\pi\circ\tilde{\phi}_i=\phi_i\circ p\qquad
\forall~i=0,\ldots,n,$$
and $\tilde{\phi}_i\!:p^{-1}(\under{w})\!\lra\!\pi^{-1}(\phi_i(\under{w}))$
is a $\Bbb{C}$-linear map for all $\under{w}\!\in\!\Bbb{C}^n$.
Thus, $\ga\!\lra\!\P$ is a {\it holomorphic rank-one vector bundle},
i.e.~a {\it holomorphic line bundle}.\\

\noindent
Each homogeneous polynomial, 
$$p=\!\!\sum_{i_0+\ldots+i_n=d}\!\!\!\!\!\!\!
a_{i_0\ldots i_n}z_0^{i_0}\ldots z_n^{i_n},$$
of degree~$d$ in $n\!+\!1$ variables determines a section $s_p$ 
of the bundle  $\ga^{*\otimes d}\!\lra\!\P$, described as follows.
At each point $\ell\!\in\!\P$, $s_p(\ell)$ is to be a map from $\ga_p$ 
to~$\Bbb{C}$ such that
$$\big\{s_p(\ell)\big\}(t\under{z})=t^d\big\{s_p(\ell)\big\}(\under{z})
\qquad\forall~\under{z}\!\in\!\ga_p=\ell.$$
Thus, we define $s_p$ by
$$\big\{s_p(\ell)\big\}\big(\ell;z_0,\ldots,z_n\big)
=p(z_0,\ldots,z_n).$$
Lemma~\ref{holom_secs_lmm} below
can be checked directly from the relevant definitions.

\begin{lmm}
\label{holom_secs_lmm}
If $p$ is a homogeneous polynomial of degree~$d$ in $n\!+\!1$ variables,
$s_p$ is a holomorphic section of the holomorphic line bundle~$\ga^{\otimes *d}$.
Conversely, if $s$ is a holomorphic section of~$\ga^{\otimes *d}$,
$s\!=\!s_p$ for some homogeneous polynomial~$p$ of degree~$d$ in $n\!+\!1$ variables.
\end{lmm}

\noindent
If $s$ is a section of a vector bundle~$V$ over a smooth manifold~$X$
and $x\!\in\!s^{-1}(0)$, the differential of~$s$ at~$x$ is 
a well-defined linear map:
$$ds\big|_x\!:T_xX\lra V_x.$$
It can be constructed using either a chart for $V$ or 
a connection in~$V$.
If $ds|_x$ is surjective, $s$~is said to be 
{\it transversal to the zero~set at~$x$}.
If $ds|_x$ is surjective for all $x\!\in\!s^{-1}(0)$,
$s$ is to be {\it transverse to the zero~set}.
If $V$ is a complex vector bundle of rank~$n$,
$X$ is a complex $n$-manifold, and 
$s$ is transversal to the zero~set at~$x\!\in\!s^{-1}(0)$,
$x$~is an isolated point of~$s^{-1}(0)$ and 
$ds|_x\!:T_xX\!\lra\!V_x$ is an $\Bbb{R}$-linear map
between complex (and thus, oriented) vector spaces.
The point~$x$ is assigned the plus sign if this map is
orientation-preserving and the minus sign otherwise.
Note that if $s$ is a holomorphic section,
$ds|_x$ is $\Bbb{C}$-linear and thus orientation-preserving.\\

\noindent 
We conclude this subsection by proving Lemma~\ref{genus_lmm}.
With notation as before,
\begin{equation}\label{genus_e1}
g\big(s^{-1}(0)\big) = \frac{2-\chi(s^{-1}(0))}{2},
\end{equation}
where $\chi(s^{-1}(0))$ is the euler characteristic of the surface $s^{-1}(0)$.
On the other hand, by Corollary~11.12 in~\cite{MiSt} and by Lemma~\ref{euler_lmm},
\begin{equation}\label{genus_e2}\begin{split}
\chi(s^{-1}(0)) &= \blr{e(Ts^{-1}(0)),s^{-1}(0)}
=\blr{c_1(T\PP)\!-\!c_1(\ga^{*\otimes d}),s^{-1}(0)}\\
&=\blr{(3a\!-\!da)\cdot da,\PP} =3d-d^2.
\end{split}\end{equation}
Lemma~\ref{genus_lmm} follows immediately from \e_ref{genus_e1} and \e_ref{genus_e2}.

\subsection{Plane Curves}
\label{curves_subs}

\noindent 
A ({\it reduced}, {\it complex}) {\it curve} ${\cal C}$ in $\PP$ is a subset of 
$\PP$ of the form
$${\cal C}={\cal C}_{\under{a}}\equiv\big\{[X,Y,Z]\!\in\!\PP\!:
\sum_{j+k+l=d}a_{jkl}X^jY^kZ^l=0\big\},$$
for some positive integer $d$ and some tuple $\under{a}=(a_{jkl})_{j+k+l=d}$
of complex numbers, not all~zero.
In other words, a curve in $\PP\!\equiv\!(\Bbb{C}^3\!-\!\{0\})/\Bbb{C}^*$
is the quotient of the zero set of a nonzero homogeneous polynomial on
$\Bbb{C}^3\!-\!\{0\}$) by the $\Bbb{C}^*$-action.
The {\it degree}~$d({\cal C})$ of the curve~${\cal C}$ in~$\PP$ is the minimal degree
of a homogeneous polynomial giving rise to~${\cal C}$.
Alternatively, $d({\cal C})$ is the positive number such that 
$$[{\cal C}]=d({\cal C})\cdot\ell\in H_2(\PP;\Bbb{Z}),$$
where $\ell$ is the homology class of a line in $\PP$.\\

\noindent
If ${\cal C}\!\subset\!\PP$ is a curve, there exists a smooth Riemann surface~$\Si$,
possibly not connected, and a holomorphic map $f\!:\Si\!\lra\!\PP$ such that
${\cal C}\!=\!f(\Si)$.
The degree of such a map~$f$ is the number $d(f)$ such that
$$f_*[\Si]=d(f)\cdot\ell\in H_2(\PP;\Bbb{Z}).$$
If ${\cal C}\!=\!f(\Si)$, $d({\cal C})\!\le\!d(f)$.
If $d({\cal C})\!=\!d(f)$, $f\!:\Si\!\lra\!{\cal C}$ is a 
{\it normalization} of~${\cal C}$.
If $f\!:\Si\!\lra\!{\cal C}$ is a normalization of~${\cal C}$,
the ({\it geometric}) {\it genus},~$g({\cal C})$, of the curve~${\cal C}$
is the genus of Riemann surface~$\Si$.\\

\noindent
The following two lemmas can be proved using basic facts from complex analysis
and algebraic topology.

\begin{lmm}
\label{normalization_lmm}
Every complex curve ${\cal C}\!\subset\!\PP$ admits a normalization 
$f\!:\Si\!\lra\!{\cal C}$.
If $f_1\!:\Si_1\!\lra\!{\cal C}$ and $f_2\!:\Si_2\!\lra\!{\cal C}$
are normalizations of ${\cal C}$, there exists a biholomorphism
$\tau\!:\!\Si_1\!\lra\!\Si_2$ such that $f_1\!=\!f_2\!\circ\!\tau$.
\end{lmm}

\begin{lmm}
\label{bezout_lmm}
If ${\cal C}_1$ and ${\cal C}_2$ are complex plane curves that intersect 
at a finite number points, then the number of intersection points counted
with appropriate positive multiplicities is $d({\cal C}_1)\cdot d({\cal C}_2)$.
\end{lmm}

\end{document}